\documentclass[12pt,tgrind,psbox,leqno]{article}
\usepackage{latexsym}
\usepackage{graphics}
\usepackage{amsmath}
\usepackage{xspace}
\usepackage{amssymb}

\def\ignore#1{{}}

\newcommand {\be}[1]{\begin{equation}\label{#1}}
\newcommand {\ee}{\end{equation}}
\newcommand {\bea}{\begin{eqnarray}}
\newcommand {\eea}{\end{eqnarray}}
\newcommand{\Pois}{\ensuremath{\operatorname{Pois}}\xspace}
\newcommand{\Be}{\ensuremath{\operatorname{Be}}\xspace}
\newcommand{\Exp}{\ensuremath{\operatorname{Exp}}\xspace}

\newcommand{\qed}{\hfill $\Box$}
\newcommand{\remark}{{\bf Remark : }}

\newcommand{\pr}{\mathbb{P}}
\newcommand{\E}{\mathbb{E}}

\newcommand{\I}{{\bf {\cal I}}}
\newcommand{\M}{{\bf {\cal M}}}

\newtheorem{theorem}{Theorem}
\newtheorem{lemma}[theorem]{Lemma}
\newtheorem{conj}{Conjecture}
\newtheorem{prop}{Proposition}
\newtheorem{coro}{Corollary}

\title{Maximum Weight Independent Sets and Matchings in Sparse Random Graphs. Exact Results using the
Local Weak Convergence Method}

\author{%
David Gamarnik \thanks{Department of Mathematical Sciences,
IBM T.J.~Watson Research Center, Yorktown Heights NY 10598, USA.
\hbox{e-mail}~{\small\texttt{gamarnik@watson.ibm.com}}}
\and
Tomasz Nowicki \thanks{Department of Mathematical Sciences,
IBM T.J.~Watson Research Center, Yorktown Heights NY 10598, USA.
\hbox{e-mail}~{\small\texttt{nowicki@watson.ibm.com}}}
\and
Grzegorz Swirscsz\thanks{
Warsaw University
\hbox{e-mail}~{\small\texttt{swirszcz@mimuw.edu.pl}}}
}

\begin{document}

\maketitle


\begin{abstract}
Let $G(n,c/n)$ and $G_r(n)$ be an  $n$-node sparse random graph and a sparse random $r$-regular graph, respectively, and
let ${\cal I}(n,r)$ and ${\cal I}(n,c)$ be the sizes of the largest independent set in $G(n,c/n)$ and $G_r(n)$.
The asymptotic value of ${\cal I}(n,c)/n$ as $n\rightarrow\infty$, can be computed using
the Karp-Sipser algorithm when $c\leq e$. For random cubic graphs, $r=3$, it is  only known
that $.432\leq\liminf_n {\cal I}(n,3)/n \leq \limsup_n {\cal I}(n,3)\leq .4591$ with high probability (w.h.p.) as $n\rightarrow\infty$,
as shown in \cite{FriezeSuen} and \cite{BollobasRegInd}, respectively.

In this paper we assume in addition that the nodes of the graph are equipped with non-negative weights,
independently generated according to some common distribution, and we consider instead the maximum weight of an independent set.
Surprisingly, we discover that for certain weight distributions,
the limit $\lim_n {\cal I}(n,c)/n$ can be computed exactly even when $c>e$, and
 $\lim_n {\cal I}(n,r)/n$ can be computed exactly for some $r\geq 2$.
For example, when the weights  are exponentially distributed with
parameter $1$, $\lim_n {\cal I}(n,2e)/n\approx .5517$, and $\lim_n
{\cal I}(n,3)/n\approx .6077$. Our results are established using
the recently developed \emph{local weak convergence} method
further reduced to a certain \emph{local optimality} property
exhibited by the models we consider. Using the developed technique
we show in addition that in the unweighted case $\liminf_n {\cal
I}(n,4)/n\geq  .3533$, which is a new lower bound. We also prove
that in any (non-random)
 graph with degree $3$ and large girth, the size of the maximum
independent set is at least $.3923n-o(n)$, improving the previous
bound $(7/18)n-o(n)$  in \cite{HopkinsStaton}. Finally, we
extend our results to maximum weight matchings in $G(n,c/n)$ and
$G_r(n)$. For the case of exponential distributions, we compute
the corresponding limits for every $c>0$ and every $r\geq 2$.
\end{abstract}

\newpage

\tableofcontents

\section{Introduction}\label{section:introduction}
Two models of random graphs considered in this paper are a sparse random graph $G(n,c/n)$ and a sparse
random regular graph $G_r(n)$. The first is a graph on $n$ nodes $\{0,1,\ldots,n-1\}\equiv [n]$, where each potential
undirected edge $(i,j), 0\leq i<j\leq n-1$ is present in the graph with probability $c/n$, independently for all
$n(n-1)/2$ edges. Here $c>0$ is a fixed constant, independent of $n$. A random $r$-regular graph
$G_r(n)$ is obtained by  fixing a constant integer $r\geq 2$  and considering a graph selected uniformly at random from the space of all
$r$-regular graphs on $n$ nodes (graphs in which every node has degree $r$). A set of nodes $V$ in a graph $G$ is defined to be
an independent set if no two nodes of $V$ are connected by an edge. Let ${\cal I}(n,c)$ and  ${\cal I}(n,r)$
denote the maximum cardinality of an independent set in $G(n,c/n)$ and $G_r(n)$ respectively.
Suppose the nodes of a graph are equipped with some non-negative weights  $W_i, 0\leq i\leq n-1$ which are
generated independently according to some common distribution $F_w(t)=\pr\{W_i\leq t\}, t\geq 0$.
Let ${\cal I}_w(n,c),{\cal I}_w(n,r)$ denote maximum weight of an independent set in $G(n,c/n)$ and
$G_r(n)$ respectively.

A  matching is a set of edges $A$ in a graph $G$ such that every node is incident to at most one edge in $A$.
Let ${\cal M}(n,c)$ and ${\cal M}(n,r)$ denote the maximum cardinality of a matching in $G(n,c/n)$ and $G_r(n)$, respectively. It
is known that  $G_r(n), r\geq 3$ has a full matching w.h.p., that is
${\cal M}(n,r)=n/2$ ($\lfloor n/2\rfloor$ for odd $n$) w.h.p. \cite{JansonBook}.
If the edges of the graph are equipped with some non-negative random weights,
then we consider instead the maximum weight of a matching  ${\cal M}_w(n,c)$ and ${\cal M}_w(n,r)$ in graphs $G(n,c/n),G_r(n)$, respectively.
The computation of ${\cal I}_w(n,c),{\cal I}_w(n,r),{\cal M}_w(n,c), {\cal M}_w(n,r)$ in the limit as $n\rightarrow\infty$
is the main subject of the present paper.

The asymptotic values of ${\cal I}(n,c)$ for $c\leq e$ and ${\cal M}(n,c)$ for all $c$ were obtained by Karp and Sipser
using a  simple greedy type algorithm in \cite{KarpSipser}. Yet it is an open problem to compute the corresponding limit for
for independent sets for the case $c>e$ or in random regular graphs or even to show that the limit exists \cite{Aldous:FavoriteProblems}.

The developments in this paper show that, surprisingly, proving the existence and the computation of the limits $\lim_n{\cal I}(n,\cdot)/n,
\lim_n{\cal M}(\cdot)/n$ is easier in the weighted case than in the unweighted case, at least for certain weight distributions. In particular,
we compute the limits for independent sets in $G_r(n), r=2,3,4$ and $G(n,c/n), c\leq 2e$, when the node weights are exponentially distributed,
and we compute the limits for matchings in $G_r(n)$ and $G(n,c/n)$ for all $r,c$, when the edge weights are exponentially distributed.
It was shown by the first  author \cite{gamarnik_LSAT}  that the limit $\lim_n{\cal M}(r,c)/n$  exists for every weight distribution
with bounded support, though the non-constructive methods employed prevented the computation of the limits.

Our method of proof is based on a powerful \emph{local weak convergence method} developed by Aldous \cite{Aldous:assignment92}, \cite{Aldous:assignment00},
Aldous and Steele \cite{AldousSteele:survey}, Steele \cite{Steele:MinSpanningTree}, further empowered by a certain
\emph{local optimality} observation derived in this paper. Local weak convergence is a  recursion
technique based on fixed points of distributional equations, which allows one to compute limits of some random combinatorial structures,
see Aldous and Bandyopadhyay \cite{AldousBandyopadhyaySurvey} for a recent survey on applications of
distributional equations and Aldous and Steele \cite{AldousSteele:survey} for a survey on the local weak convergence
method. In particular, the method is used to compute maximum weight matching on  a random tree, when the weights are exponentially distributed.
The tree structure was essential in \cite{AldousSteele:survey} for certain computations and the approach
does not extend directly to graphs like $G(n,c/n)$ with $c>1$, where the convenience of a tree structure is lost
due to the presence of a giant component. It was conjectured in \cite{AldousSteele:survey} that a some long-range independence property might
be helpful to deal with this difficulty. The present paper partially answers this qualitative conjecture in a positive way.
We introduce a certain operator $T$ acting on the space of distribution functions. We prove a certain  local optimality property
stating that, for example, for independent sets, whether a given node $i$ belongs to the maximum weight independent set is asymptotically
independent from the portion of the graph outside  a constant size neighborhood of $i$, iff $T^2$ has a unique fixed point
distribution. Moreover, when $T^2$ does have the unique fixed point,
the size of the extremal object (say maximum weight independent set) can be
derived from a fixed point of an operator $T$. The computations of fixed points is tedious, but simple in principle and the groundwork for
that was already done in \cite{AldousSteele:survey}. We hope that the long-range independence holds  in
other random combinatorial structures as well.

The issue of long-range independence of random combinatorial objects is addressed in a somewhat different, statistical physics
context in Mossel \cite{EMosselSurvey}, Brightwell and Winkler \cite{BrightwellWinkler}, Rozikov and Suhov \cite{RozikovSuhov}, Martin \cite{Martin},
where independent sets (hard-core model) are considered on infinite regular trees, weighted
by the Gibbs measure. The long-range interaction between nodes at a large distance
is investigated with respect to this measure using  the notion of reconstruction.
It would be interesting to investigate the connections between the two models.

In a different setting Talagrand \cite{talagrandRASSGN} proves a certain long-range independence property for the random assignment problem, where
the usual min-weight matching is replaced by a partition function on the space of feasible matchings.
He uses a rigorous mathematical version of the cavity method, which
originated in physics, to prove that the spins (edges of the matching) are
asymptotically independent as the size of the problem increases. The particular form of the long-range independence is similar to the one we obtain, refer
to Theorem \ref{theorem:Gaussian} below, and in fact the cavity method, which is based on "knocking" out certain spins from the system and analyzing
the relative change of the size of the extremal object, has some similarity with the local weak convergence method, which is also based on considering
extremal objects (say independent sets) with one or several nodes excluded. It seems worth investigating whether there is  a formal connection
between the two methods. Finally, we  refer the reader to Hartmann and Weigt \cite{HartmannWeigt} who derive the same result as
Karp and Sipser for independent sets using non-rigorous arguments from statistical physics.

The rest of the paper is organized as follows. In the following
section we describe some prior results on maximum
independent sets and matchings in random graphs. Our main theorems
are given in Section \ref{section:results}. The operator $T$,
fixed point equations  and long-range independence issues are
discussed in Section \ref{section:FixedPoints}. The main results
are proven in Section \ref{section:applications to G}. Some conclusions are
in Section \ref{section:conclusions}

We finish this section with some notational conventions. $\Exp(\mu),\Pois(\lambda),\Be(p)$
denote respectively exponential, Poisson and  Bernoulli  distributions with parameters
$\mu,\lambda>0, 0\leq p\leq 1$.

\section{Prior work and open questions}\label{section:background}
It is known and simple to prove that $\I(n,r)=\I(n,c)=\Theta(n)$ w.h.p. for any constants $r\geq 1,c>0$.
Moreover, it is known that, w.h.p.,
$6\log(3/2)-2=.432\ldots\leq \liminf_n\I(n,3)/n\leq \limsup_n\I(n,3)/n\leq .4591$. The lower bound is due to Frieze and Suen \cite{FriezeSuen},
and the upper bound is due to Bollobas \cite{BollobasRegInd}. The upper bound is generalized for any $r\geq 0$
and uses a very ingenious construction of random regular graphs via matching and random grouping, \cite{BollobasRandReg},
\cite{JansonBook}.

It is natural to expect that the following is true, which unfortunately remains only a conjecture,
appearing in several places, most recently in \cite{Aldous:FavoriteProblems} and \cite{AldousSteele:survey}.

\begin{conj}\label{conj:indset}
For every $c>0$ and $r\geq 3$ the limits
\[
\lim_{n\rightarrow\infty}{\E[\I(n,c)]\over n}, \qquad \lim_{n\rightarrow\infty}{\E[\I(n,r)]\over n}
\]
exist.
\end{conj}
The existence of these limits also implies the convergence to the same limits w.h.p. by applying Azuma's inequality,
see \cite{JansonBook} for the statement and the applicability of this inequality.

The limit $\lim_n\E\I(n,c)/n$ is known to exist for $c\leq e$ and can be computed using the Karp-Sipser  \cite{KarpSipser} algorithm  for
maximum matching $\M(n,c)$. We describe the algorithm first in high level terms, and then state the result and its implications to
maximum independent sets. The algorithm proceeds in two stages. In the first stage any leaf $v$ in the graph $G(n,c/n)$ is selected.
The edge incident to this leaf is selected into a matching, and all the other edges incident to the parent of $v$ are deleted. This is
repeated until no leaves are left. In the second stage simply the  largest matching $M^*$ in the remaining graph $G^*\subset G$ is selected
and added to the matching constructed in first stage. It is a simple exercise to prove that the resulting matching is optimal.
Karp and Sipser prove the following result about the size $\M(KS)$ of the obtained matching.

\begin{theorem}[Karp and Sipser \cite{KarpSipser}]\label{theorem:KarpSipser}
The constructed matching is optimal, $\M(KS)=\M(n,c)$ with probability one. Moreover,
\begin{enumerate}
\item The maximum matching
satisfies
\be{eq:MatchingCard}
 \lim_{n\rightarrow\infty}{\E[\M(n,c)]\over n}=1-{\gamma^*(c)+\gamma^{**}(c)+c\gamma^*(c)\gamma^{**}(c)\over 2},
\ee
where $\gamma^*(c)$ is the smallest solution of the equation $x=\exp(-c\exp(-cx))$ and $\gamma^{**}(c)=\exp(-c\gamma^*(c))$.
\item When $c\leq e$, the equation $x=\exp(-c\exp(-cx))$ has a unique solution $\gamma(c)$, and
$|G^*|=o(n)$. That is the matching constructed in the first stage is asymptotically optimal as $n\rightarrow\infty$.
\end{enumerate}
\end{theorem}

The different behavior for $c\leq e$ and $c>e$ is called $e$-cutoff phenomena.
Theorem \ref{theorem:KarpSipser} was strengthened later by by Aronson, Frieze and Pittel \cite{AronsonFriezePittel},
who obtain bounds on convergence (\ref{eq:MatchingCard}).
It is the second part of the theorem above that can be used for the analysis of $\I(n,c)$. Observe, that in the
first stage of the Karp-Sipser algorithm for every selected leaf if one takes a parent of this leaf instead of an
edge between the leaf and its parent, one obtains a node set which is an edge-cover set in the graph $G\setminus G^*$
constructed in the first stage. That is every edge in $G\setminus G^*$ is incident to at least one node in this cover.
It is a simple exercise to prove that this is in fact a minimum node cover in $G\setminus G^*$. Its complement is a
maximum independent set in $G\setminus G^*$. Since $|G^*|=o(n)$ when $c\leq e$, one obtains then the following result.

\begin{coro}\label{corollary:IndSetc<e}
When $c\leq e$:
\be{eq:IndSetCard_c<e}
\lim_{n\rightarrow\infty}{\E[\I(n,c)]\over n}={2\gamma(c)+\gamma^2(c)\over 2}.
\ee
\end{coro}

The Karp-Sipser algorithm hinges strongly on working with leaves and thus is not applicable to random
regular graphs. Moreover, if the edges or the nodes of the graph $G(n,c/n)$ are equipped with weights then
the Karp-Sipser  algorithm clearly can produce a strictly suboptimal solution and cannot be used in
our setting of weighted nodes and edges. Also, when the edges of
$G_r(n)$ are equipped with weights, the problem of computing maximum weight matching becomes non-trivial,
as opposed to the unweighted case when the full matching exists w.h.p.

In a somewhat different domain of extremal combinatorics the following result was established by Hopkins and
Staton \cite{HopkinsStaton}. A girth of a graph is the size of the smallest cycle. It is shown in \cite{HopkinsStaton} that the size
of a largest independent set in an $n$-node graph with largest degree $3$ and large girth is asymptotically at least $(7/18)n-o(n)$.
The techniques we employ in this  paper allow us to improved this lower bound.

\section{Main results}\label{section:results}
We begin by introducing the key technique for our analysis -- recursive distributional equations and its fixed point solutions.
This technique was introduced by Aldous \cite{Aldous:assignment92}, \cite{Aldous:assignment00} in the context of
solving $\zeta(2)$ limit conjecture for random minimal assignment problem,
and was further developed in  Aldous and Steele \cite{AldousSteele:survey}, Steele \cite{Steele:MinSpanningTree}, Aldous and
Bandyopadhyay \cite{AldousBandyopadhyaySurvey}, Gamarnik \cite{gamarnik_LSAT}. Let  $W$ be a non-negative random variable with a distribution function $F_w(t)=\pr(W\leq t)$.
We consider four operators $T=T_{{\cal I},r},T_{{\cal I},c},T_{{\cal M},r},T_{{\cal M},c}$ acting on
the space of distribution functions $F(t),t\geq 0$, where $c>0$ is a fixed constant and $r\geq 2$ is a fixed integer.

\begin{enumerate}
\item Given $W$  distributed according to $F_w$ (we write simply $W \sim F_w$), and given a distribution function $F=F(t)$,
let $B_1,B_2,\ldots,B_r \sim F$
be generated independently. Then $T_{{\cal I},r}:F\rightarrow F'$, where $F'$ is the distribution function
of $B'$ defined by
\be{eq:RecursionIndr}
B'=\max(0,W-\sum_{1\leq i\leq r}B_i).
\ee

\item Under the same setting as above, let $B_1,\ldots,B_m \sim F$, where $m$ is a random variable distributed according
to a Poisson distribution with parameter $c$, independently from $W,B_i$.
Then $T_{{\cal I},c}:F\rightarrow F'$, where $F'$ is the distribution function
of $B'$ defined by
\be{eq:RecursionIndc}
B'=\max(0,W-\sum_{1\leq i\leq m}B_i),
\ee
when $m\geq 1$ and $B'=W$ when $m=0$. For simplicity we will identify the sum above
with zero when $m=0$.
\item Let  $W_1,\ldots,W_r\sim F_w$, $B_1,\ldots,B_r\sim F$.
Then $T_{{\cal M},r}:F\rightarrow F'$, where $F'$ is the distribution function
of $B'$ defined by
\be{eq:RecursionMr}
B'=\max_{1\leq i\leq r}(0,W_i-B_i).
\ee

\item Finally, let $W_1,\ldots,W_m\sim F_w$, $B_1,\ldots,B_m\sim F$, where $m\sim\Pois(c)$, independent
from $W_i,B_i$. Then $T_{{\cal M},c}:F\rightarrow F'$, where $F'$ is the distribution function
of $B'$ defined by
\be{eq:RecursionMc}
B'=\max_{1\leq i\leq m}(0, W_i-B_i),
\ee
when $m\geq 1$ and $B'=0$ when $m=0$. Again, for simplicity, we assume that  $\max$ expression above is zero when $m=0$.
\end{enumerate}

A distribution function $F$ is defined to be a fixed point distribution of an operator $T$ if $T(F)=F$.

We now state the main result of this paper. Recall, that a distribution function $F(t)$ is defined to be continuous
(atom free) if for every $x$ in its support $\lim_{\epsilon\rightarrow 0}(F(x+\epsilon)-F(x-\epsilon))=0$. Equivalently,
for $B\sim F$ and every $x$, $\pr(B=x)=0$. We use $1\{\cdot\}$ to denote the indicator function.

\begin{theorem}\label{theorem:main}
Let $F_w$ be a continuous non-negative distribution function.  For $r\geq 1$ if the operator $T^2_{{\cal I},r-1}$
has a unique fixed point distribution function $F^*$,
then, w.h.p.
\be{eq:mainIndr}
\lim_n {{\cal I}_w(n,r)\over n}=\E[W\,1\{W-\sum_{1\leq i\leq r}B_i>0\}],
\ee
where $W\sim F_w$,  $B_i\sim F^*$, and $W,B_i$ are independent.
When $G_r(n)$ is replaced by $G(n,c/n)$, the same result holds for $T=T_{{\cal I},c}$,
except the sum in the right-hand side of (\ref{eq:mainIndr}) is $\sum_{1\leq i\leq m}B_i$ and $m\sim \Pois(c)$.

Finally, the similar results hold for ${\cal M}_w(n,r)$ and ${\cal M}_w(n,c)$ in
$G_r(n)$ and $G(n,c)$, for $T=T_{{\cal M},r-1}$ and $T=T_{{\cal M},c}$, respectively, whenever the corresponding
operator $T$ is such that $T^2$ has the unique fixed point distribution $F^*$. The corresponding limits are
\be{eq:mainMr}
\lim_n {{\cal M}_w(n,r)\over n}={1\over 2}\E[\sum_{1\leq i\leq r}W_i\,1\{W_i-B_i=\max_{1\leq j\leq r}(W_j-B_j)>0\}],
\ee
where $W_i\sim F_w$, $B_i\sim F^*$,
and
\be{eq:mainMc}
\lim_n {{\cal M}_w(n,c)\over n}={1\over 2}\E[\sum_{i\leq m}W_i\,1\{W_i-B_i=\max_{j\leq m}(W_j-B_j)>0\}],
\ee
where $W_i\sim F_w, B_i\sim F^*, m\sim \Pois(c)$.
\end{theorem}

For  $G=G_r(n),G(n,c/n)$ with the weights on nodes given by a distribution function $F_w$,
let ${\cal IN}_w(n,r),$ ${\cal IN}_w(n,c)$ denote the cardinality of the maximum weight independent
set in $G$ which achieves the maximum weight. In case $F_w$ is continuous, the maximum weight independent
set is uniquely defined, so ${\cal IN}_w(n,r),{\cal IN}_w(n,c)$ are well-defined as well.
Clearly, ${\cal I}(n,r)\geq {\cal IN}_w(n,r), {\cal I}(n,c)\geq {\cal IN}_w(n,c)$. ${\cal MN}_w(n,r)$
and ${\cal MN}_w(n,c)$ are defined similarly. When $T^2_{{\cal I},r}$ or $T^2_{{\cal I},r}$ have the unique
fixed point it is also possible to compute asymptotically ${\cal IN}_w(n,r), {\cal IN}_w(n,c),$ ${\cal MN}_w(n,r), {\cal MN}_w(n,c)$.

\begin{coro}\label{coro:main}
Under the setting of Theorem \ref{theorem:main}, w.h.p.
\be{eq:mainIndrCard}
\lim_n {{\cal IN}_w(n,r)\over n}=\E[1\{W-\sum_{1\leq i\leq r}B_i>0\}],
\ee
where $W\sim F_w$,  $B_i\sim F^*$, $W,B_i$ independent, and $F^*$ is the unique fixed point distribution of $T^2=T^2_{{\cal I},r-1}$,
and
\be{eq:mainMrCard}
\lim_n {{\cal MN}_w(n,r)\over n}={1\over 2}\E[1\{\max_{1\leq j\leq r}(W_j-B_j)>0\}],
\ee
where $W_i,\sim F_w$,  $B_i\sim F^*$, $W_i,B_i$ independent, and $F^*$ is the unique fixed point distribution of $T^2=T^2_{{\cal M},r-1}$
and $W,B_i$ are independent,

Similar results hold for $T=T_{{\cal I},c},T_{{\cal M},c}$ where again in the sum $\sum_{1\leq j\leq r}$
we substitute $r$ with random $m\sim \Pois(c)$.
\end{coro}

The Theorem \ref{theorem:main}  is the core result of this paper. It will allow us to obtain several
interesting corollaries, which we state below.

\begin{theorem}\label{theorem:mainEXPr}
Suppose the weights of the nodes and edges of the graphs $G=G_r(n)$ and $G=G(n,c)$ are distributed as $\Exp(1)$.
Then
\begin{enumerate}
\item $T^2_{{\cal I},r-1}$ has a unique fixed point distribution $F^*$ iff $r\leq 4$.
In this case, w.h.p.
\be{eq:mainIndr=3}
\lim_n {{\cal I}_w(n,r)\over n}={(1-b)(r-r b+2b+2)\over 4},
\ee
where $b$ is the unique solution of $b=1-({1+b\over 2})^{r-1}$.
In particular, w.h.p.
\be{eq:mainIndr=2,3,4}
\lim_n {{\cal I}_w(n,2)\over n}={2\over 3},\qquad \lim_n {{\cal I}_w(n,3)\over n}\approx .6077,\qquad \lim_n {{\cal I}_w(n,4)\over n}\approx .4974,
\ee
\item $T^2_{{\cal I},c}$ has a unique fixed point distribution $F^*$ iff $c\leq 2e$.
In this case, w.h.p.
\be{eq:mainIndc}
\lim_n {{\cal I}_w(n,c)\over n}=(1-b)(1+{c(1-b)\over 4}),
\ee
where $b$ is the unique solution of $1-b=e^{-{c\over 2}(1-b)}$.
In particular, when $c=2e$, this limit is $\approx .5517$.
\item $T^2_{{\cal M},r-1}$ has a unique fixed point $F^*$ for every  $r\geq 2$.
Moreover, w.h.p.
\be{eq:mainMr=3}
\lim_n {{\cal M}_w(n,r)\over n}=r(b^{r-1}+1)\int_0^{\infty}te^{-t}(1-e^{-t}(1-b)))^{r-1}dt-r\int_0^{\infty}te^{-t}(1-e^{-t}(1-b)))^{2r-2}dt,
\ee
where $b$ is the unique solution of $b=1-{1-b^r\over r(1-b)}$.

\item $T^2_{{\cal M},c}$ has a unique fixed point $F^*$ for all $c>0$. Moreover, for every $c>0$ w.h.p.
\be{eq:mainMcc}
\lim_n {{\cal M}_w(n,c)\over n}={c\over 2}(e^{cb-c}+1)\int_0^{\infty}te^{-t-c(1-b)e^{-t}}dt-
{c\over 2}\int_0^{\infty}te^{-t-2c(1-b)e^{-t}+c(1-b)^2e^{-2t}}dt,
\ee
where $b$ is the unique solution of $1-e^{-cb}=c(1-b)^2$.
\end{enumerate}
\end{theorem}

The result above generalizes to the case when  $\Exp(1)$ is replaced by $\Exp(\mu)$ for any
$\mu>0$, since for any $\alpha>0$ and $W\sim \Exp(1)$, $\alpha W\sim \Exp(1/\alpha)$. The expression in
(\ref{eq:mainMr=3}) involving integrals is similar to the one found in \cite{AldousSteele:survey} for
maximum weight matching on a tree. It is a pleasant surprise, though, that the answers for
independent sets are derived in closed form.

Part 2 of  Theorem \ref{theorem:mainEXPr} leads to an interesting phase transition behavior. For $F_w=\Exp(1)$
our result says that the value of $c=2e$ is a phase transition point for the operator $T^2_{{\cal I},c}$
where for $c\leq 2e$ the operator has a unique fixed point distribution, but for $c>2e$ the fixed point distribution
is not unique. Contrast this with $e$-cutoff phenomena described above. It turns out (see Theorems
\ref{theorem:long range independence}, \ref{theorem:long range dependence} below) that this phase transition is directly
related to some long range independence/dependence property of the maximum weight independent sets in the underlying graph
$G(n,c/n)$. Curiously, no such phase transition occurs for maximum weight matchings.

As a sanity check, let us show directly the validity of (\ref{eq:mainIndr=3}) for the case $r=1$. The answer
given by the formula is $3/4$. It is, though, easy
to compute ${\cal I}_w(n,1)$ exactly. The graph is a collection of $n/2$ isolated edges. For each edge we simply select the
the incident node with larger weight, which by memoryless property of $\Exp$ is $1/2+1=3/2$. The limit ${\cal I}_w(n,1)/n\rightarrow 3/4$ then checks.

Note that $G_2(n)$ can be represented as a collection
of disjoint cycles on $n$ elements. W.h.p. these cycles have decreasing lengths,
starting from $\Theta(n^{1\over 2})$. Theorems \ref{theorem:main} and \ref{theorem:mainEXPr} also hold also for a simpler
model of a $2$-regular graph -- $n$-cycle. The $n$ nodes $0,1,\ldots,n-1$ and edges $(0,1),\ldots,(n-2,n-1),(n-1,0)$
of this cycle are assumed to have weights distributed according to some distribution function $F_w$.
Let $I_w(n,{\rm cycle})$ and $M_w(n,{\rm cycle})$ denote respectively the maximum weight of an independent set and
the maximum weight of a matching.

\begin{coro}\label{coro:cycle}
Suppose $T^2_{{\cal I},2}$ has a unique fixed point. Then (\ref{eq:mainIndr}), (\ref{eq:mainMr}) hold when $r=2$
and ${\cal I}_w(n,{\rm cycle})$ replaces ${\cal I}_w(n,2)$ and ${\cal M}_w(n,{\rm cycle})$ replaces ${\cal M}_w(n,2)$.
When $F_w=\Exp(1)$, w.h.p.
\be{eq:cycle}
\lim_n {{\cal I}_w(n,{\rm cycle})\over n}=\lim_n {{\cal M}_w(n,{\rm cycle})\over n}={2\over 3}.
\ee
\end{coro}

We can also compute the cardinality of independent sets which achieve the maximum weight, when $F_w=\Exp(1)$.
\begin{coro}\label{coro:mainCardinality}
When $F_w=\Exp(1)$, w.h.p.
\be{eq:mainIndrCard=2,3,4}
\liminf_n {{\cal I}(n,2)\over n}\geq \lim_n {{\cal IN}_w(n,2)\over n}={4\over 9},
\qquad \liminf_n {{\cal I}(n,3)\over n}\geq\lim_n {{\cal IN}_w(n,3)\over n}\approx .3923,
\ee
\[
\qquad \liminf_n {{\cal I}(n,4)\over n}\geq\lim_n {{\cal IN}_w(n,4)\over n}\approx .3533.
\]
\end{coro}
The last part of the corollary above implies  a new lower bound on the size of the largest independent
sets in a random $4$-regular graph.

\begin{coro}\label{coro:4regular}
The cardinality of the largest independent set in $G_4(n)$ is at least $.3533n-o(n)$.
\end{coro}
We note, however, that for $r=3$ our lower bound $.3923n-o(n)$
is weaker than $.432n-o(n)$ established in \cite{FriezeSuen}.
Yet, it is still very useful
since the approach allows us the improve the following result of Hopkins and Staton \cite{HopkinsStaton}.
Let ${\cal G}(n,r,d)$ denote the class of all (non-random) graphs on $n$ nodes, with maximum degree $r$ and girth at least $d$.
For any $G\in {\cal G}(n,r,d)$ let ${\cal I}(G)$ denote the size of the largest independent set in $G$.
Hopkins and Staton proved that $\liminf_{n,d}\min_{G\in {\cal G}(n,3,d)}{\cal I}(G)/n\geq 7/18\approx .3887$.
Our techniques allow us to obtain the following improvement.

\begin{theorem}\label{theorem:WorstcaseGirth}
\[
\liminf_{n,d}\min_{G\in {\cal G}(n,3,d)}{{\cal I}(G)\over n}\geq .3923.
\]
\end{theorem}

An important implication of the uniqueness of the fixed point distribution of
$T^2$ for the types of $T$ described above, is that the uniqueness implies a certain long-range independence
 properties of the  structures we consider.
The following theorem makes this notion precise. While the theorem is not used directly in this paper, we believe it is
interesting by itself. Below $G$ is again one of the graphs $G_r(n)$ or $G(n,c)$.
Let $E$ denote the (random) edge set of $G$.

\begin{theorem}\label{theorem:Gaussian}
Let $T$ be one of the four operators (\ref{eq:RecursionIndr}),(\ref{eq:RecursionIndc}),(\ref{eq:RecursionMr}),(\ref{eq:RecursionMc})
with respect to some continuous distribution function $F_w$,
and let ${\cal C}_w(n)={\cal I}_w(n,r),{\cal I}_w(n,c),{\cal M}_w(n,r)$
or ${\cal M}_w(n,c)$. Denote by ${\cal O}(n)$  the subset of $[n]$ or $E$ (depending on a context), which achieves
${\cal C}_w(n)$. Select two elements $i,j$ of $[n]$ or $E$ uniformly at random. If $T^2$ has  a unique fixed point distribution, then
\be{eq:LongRangeIndependence}
\pr(i,j\in{\cal O}(n))\rightarrow \pr(i\in{\cal O}(n))\pr(j\in{\cal O}(n)),
\ee
as $n\rightarrow\infty$.
\end{theorem}

Recall that for each of the four objects $|{\cal C}(n)|=\Theta(n)$ w.h.p. As a result the values
$\pr(v\in{\cal O}(n))$ do not vanish and the first part of the theorem above does have a non-trivial
content. In fact we will show a much  stronger result stating that, when $T^2$ has a unique fixed point,
the event $v\in{\cal O}(n)$ is almost independent of the entire graph $G$ outside a depth-$d$ graph-theoretic
neighborhood of $v$, when $d$ is a sufficiently large constant integer.

\section{Fixed points of the operator $T^2$ and the long-range independence}\label{section:FixedPoints}

\subsection{Maximum weight independent sets and matchings in trees. Fixed points of $T^2$ and the bonus function}\label{subsection:Trees}
We start by analyzing operator $T$ -- one of the four operators introduced in the previous section.
Given two distribution functions $F_1,F_2$ defined on $[0,\infty)$, we  say  that $F_2$ stochastically
dominates $F_1$ and write $F_1\prec F_2$ if $F_1(t)\geq F_2(t)$ for every $t\geq 0$.
A sequence of distribution functions $F_n$ is defined to converge weakly to a distribution function
$F$ (written $F_n\Rightarrow F$) if $\lim_n F_n(t)=F(t)$ for every $t$ which is a point of continuity of $F$.

\begin{lemma}\label{lemma:continuity}
The operators $T=T_{{\cal I},r},T_{{\cal I},c},T_{{\cal M},r},T_{{\cal I},c}$ are continuous
with respect to the weak convergence. That is, given a sequence of distributions $F,F_s, s=1,2,\ldots\,\,$, if $F_s\Rightarrow F$
then  $T(F_s)\Rightarrow T(F)$.
\end{lemma}

\begin{proof}
The proof is almost immediate for $T=T_{{\cal I},r}$. Since the summation, subtraction and $\max$ operations are
continuous functions then the assertion holds. Here we use the fact that  for any continuous function $f$, $F_n\Rightarrow F$
implies $f(F_n)\Rightarrow f(F)$ (see Continuous Mapping Theorem in \cite{durrett}).
The proof for the case $T=T_{{\cal I},c}$ is slightly more
subtle since we are dealing with sum of randomly many elements. Let $X^{(n)}_1,X^{(n)}_2,\ldots,X^{(n)}_m$
be distributed independently according to $F_n$ and let $m$ be random variable with $\Pois(c)$ distribution.
Let $X_1,\ldots,X_m$ be distributed according to $F$.
Define $X^{(n)}=\max(0,W-\sum_{1\leq i\leq m}X^{(n)}_i)$ and $X=\max(0,W-\sum_{1\leq i\leq m}X_r)$.
Fix arbitrary $\epsilon>0$ and $m_0$ such that
$m$ does not exceed $m_0$ with probability at least $1-\epsilon$. For any fixed $t\geq 0$
\[
|\pr(X^{(n)}\leq t)-\pr(X^{(n)}\leq t|m\leq m_0)\pr(m\leq m_0)|\leq \epsilon,
\]
and
\[
|\pr(X\leq t)-\pr(X\leq t|m\leq m_0)\pr(m\leq m_0)|\leq \epsilon,
\]
Now, we have $|\pr(X^{(n)}\leq t|m\leq m_0)-\pr(X\leq t|m\leq m_0)|\leq \epsilon$ for
sufficiently large $n$ since $F_n\Rightarrow F$ and $m$ is conditioned to be at most $m_0$.
Combining, we obtain $|\pr(X^{(n)}\leq t)-\pr(X\leq t)|\leq 3\epsilon$ for sufficiently large $n$.
This completes the proof for $T=T_{{\cal I},c}$. The proofs for $T=T_{{\cal M},r},T=T_{{\cal M},c}$
are similar. \qed
\end{proof}

Let $0$ denote (for simplicity) the distribution function of a random variable $X$ which is zero w.p.1.
Let $W_r=\max_{1\leq i\leq r}W_i$ where $W_i\sim F_w$ are independent. Let also
$W_c=\max_{1\leq i\leq m}W_i$, where $W_i\sim F_w$ are independent and $m\sim \Pois(c)$.
Denote by $F_{w,r}$ and $F_{w,c}$ the distribution functions of $W_r$ and $W_c$, respectively.
\begin{prop}\label{prop:Tconvergence}
Fix $r\geq 1,c>0$,  a distribution function $F_w(t),t\geq 0$ and $T=T_{{\cal I},r}$ or $T_{{\cal I},c}$.
As $s\rightarrow\infty$, the two sequences of distributions
$T^{2s}(0),T^{2s}(F_w)$ weakly converge respectively to some distribution functions $F_{**},F^{**}$ which are fixed points
of the operator $T^2$. For any distribution function $F_0=F_0(t), t\geq 0$, $T^{2s}(0)\prec T^{2s}(F_0)\prec T^{2s}(F_w)$
and $T^{2s+1}(F_w)\prec T^{2s+1}(F_0)\prec T^{2s+1}(0)$ for all $s=1,2,\ldots$.

If the operator $T$
is such that $T^2$ has a unique fixed point (and  $F_{**}=F^{**}\equiv F^*$), then for
any distribution function $F_0=F_0(t), t\geq 0$, $T^s(F_0), s=1,2,\ldots$ converges to $F^*$ as $s\rightarrow\infty$.
In particular, $T^s(0), T^s(F_w)\rightarrow F^*$. Moreover, $F^*$ is also the unique fixed point of $T$.

When $T=T_{{\cal M},r}$ or $T=T_{{\cal M},c}$ the same result holds with $F_{w,r}$ and $F_{w,c}$ respectively
replacing $F_w$.

\end{prop}

\begin{proof}
Let $F=F(t),t\geq 0$ be any distribution function corresponding to any non-negative random variable.
It follows immediately from the definitions that
\be{eq:F}
T_{{\cal I},r}(F)\prec F_w,\qquad
T_{{\cal I},c}(F)\prec F_w,\qquad
T_{{\cal M},r}(F)\prec F_{w,r},\qquad
T_{{\cal M},c}(F)\prec F_{w,c}.
\ee
Observe that the each of the four operators $T$ above is anti-monotone. That is if $F_1\prec F_2$ for some distribution functions
$F_1,F_2$, then $T(F_1)\succ T(F_2)$. Applying this twice we obtain $T^2(F_1)\prec T^2(F_2)$.
Then $0\prec T^2(0)\prec T^4(0)\prec\cdots$, and this sequence weakly converges to
some function $F_{**}=F(t)_{**}\leq 1$, since the values of $T^{2s}(0)$ at each fixed $t$ are decreasing and are
bounded below by $0$.  Note that $F_{**}(t)$ is a non-decreasing function of $t$,
since this is the case for each distribution $T^{2s}(0)$. Finally, from (\ref{eq:F}) we have
that $\lim_{t\rightarrow\infty}F_{**}(t)=1$.
Thus $F_{**}(t)$ is a distribution function. As $T^{2s}(0)\Rightarrow F_{**}$ and by Lemma \ref{lemma:continuity}, $T$ is continuous,
we have $T^2(F_{**})=F_{**}$.

We now fix $T=T_{{\cal I},r}$. The proof of the other three cases is very similar.
From (\ref{eq:F}), by taking  $F=T_{{\cal I},r}(F_w)$ we obtain $T^2(F_w)\prec F_w$.
Then, by monotonicity of $T^2$ we obtain $F_w\succ T^2(F_w)\succ\cdots\succ T^{2s}(F_w)$
and the sequence $T^{2s}(F_w)$ converges weakly to some  function $F^{**}$. We repeat the arguments
above to show that $F^{**}$ is actually a distribution function.
Again applying Lemma \ref{lemma:continuity},
we conclude $T^2(F^{**})=F^{**}$.

Suppose now $T^2$ has a unique fixed point $F^*=F^{**}=F_{**}$.
For any distribution function $F_0$ we have $0\prec T(F_0),T^2(F_0)\prec F_w$, where again we use
(\ref{eq:F}). Applying the monotonicity $T^{2s}(0)\prec T^{2s+1}(F_0),T^{2s+2}(F_0)\prec T^{2s}(F_w)$.
But since $T^{2s}(0),T^{2s}(F_w)\Rightarrow F^*$, then $T^s(F_0)\Rightarrow F^*$. In particular, by taking $F_0=0$
and $F_0=F_w$, we obtain $T^s(0),T^s(F_w)\rightarrow F^*$.
Finally, we obtain $T^{2s}(T(F^*))\rightarrow F^*$. But $T^{2s}(T(F^*))=T(T^{2s}(F^*))=T(F^*)$. Thus $T(F^*)=F^*$. Clearly
it is the unique fixed point of $T$ since any fixed point of $T$ is also a fixed point of $T^2$. \qed
\end{proof}

We now switch to analyzing the maximum weight independent set problem on a tree. The derivation
here repeats the development in \cite{AldousSteele:survey} for maximum weight matching in random trees.
We highlight important differences where appropriate.

Suppose we have a (non-random) finite tree $H$ with nodes $0,1,\ldots,h=|H|-1$, with a fixed root $0$. The nodes of this tree are
equipped with some (non-random) weights $W_0,W_1,\ldots,W_h\geq 0$. For any node $i\in H$, let $H(i)$ denote the subtree rooted at $i$ consisting of all the
descendants of $i$. In particular, $H(0)=H$. Let ${\cal I}_{H(i)}$ denote the maximum weight of an independent set in  $H(i)$
and let $B_{H(i)}={\cal I}_{H(i)}-\sum_j{\cal I}_{H(j)}$, where the sum runs over nodes $j$ which are children of $i$.
If $i$ has no children then $B_{H(i)}$ is simply ${\cal I}_{H(i)}=W_i$.
Observe, that $B_{H(i)}$ is also a difference between ${\cal I}_{H(i)}$ and the maximum weight of an independent
set in $H(i)$, which is not allowed to use node $i$. Clearly, $0\leq B_{H(i)}\leq W_i$. The value
$B_{H(i)}$ was considered in \cite{AldousSteele:survey} in the context of maximum weight matchings and was
referred to as a \emph{bonus} of a node $i$ in tree $H(i)$.
W.l.g. denote by $1,\ldots,m$ the children of the root node $0$.

\begin{lemma}\label{lemma:recursionIndSet}
\be{eq:recursionIndSet}
B_{H(0)}=\max(0,W_0-\sum_{1\leq i\leq m}B_{H(i)}).
\ee
Moreover, if $W_0>\sum_{1\leq i\leq m}B_{H(i)}$ (that is if $B_{H(0)}>0$) then the maximum weight independent
set  must contain node $0$. If $W_0<\sum_{1\leq i\leq m}B_{H(i)}$ then the maximum weight independent
set  does not contain the node $0$.
\end{lemma}

\remark There might be several independent sets in $H$ which achieve ${\cal I}_H$. The second part of the
lemma refers to \emph{any} independent set achieving maximum weight. Also, the statement of the lemma
applies as well to every  node $i$ with respect to its tree $H(i)$. That is let $j_1,\ldots,j_l$
be the children of a node $i$. Then the lemma claims $B_{H(i)}=\max(0,W_0-\sum_{1\leq i\leq l}B_{H(j_i)}).$

\begin{proof}
Consider an independent set $V\subset H$ which achieves the maximum weight. We take an arbitrary
such in case there are many. If $0\in V$ then $i\notin V$ for $1\leq i\leq m$. Then $V$ is obtained
by taking maximum independent sets $V_i$ in $H(i)$ such that $i\notin V(i)$. By definition the
weight of such $V_i$ is ${\cal I}_{H(i)}-B_{H(i)}$. On the other hand, if $0\notin V$, then a
maximum weight of an independent set in $H$ is obtained simply as $\sum_{1\leq i\leq m}{\cal I}_{H(i)}$.
We conclude that
\be{eq:recursion1}
{\cal I}_H=\max(W_0+\sum_{1\leq i\leq m}({\cal I}_{H(i)}-B_{H(i)}),\sum_{1\leq i\leq m}{\cal I}_{H(i)}).
\ee
Recall that $B_H={\cal I}_H-\sum_{1\leq i\leq m}{\cal I}_{H(i)}$. Subtracting $\sum_{1\leq i\leq m}{\cal I}_{H(i)}$
from both sides of (\ref{eq:recursion1}) we obtain (\ref{eq:recursionIndSet}). The second part of the lemma
follows directly from the discussion above. \qed
\end{proof}

A similar development is possible for maximum weight matching. Suppose the edges of the tree $H$
are equipped with weights $W_{i,j}$. Let ${\cal M}_{H(i)}$ denote the maximum weight of a
matching in $H(i)$, and let $B_{H(i)}$ denote the difference
between ${\cal M}_{H(i)}$ and the maximum weight of a matching in $H(i)$ which is not allowed
to include any edge in $H(i)$ incident to $i$. Again $1,2,\ldots,m$ are assumed to be the children of the root $0$.

\begin{lemma}\label{lemma:recursionM}
\be{eq:recursionM}
B_H=\max(0,\max_{1\leq i\leq m}(W_{0,i}-B_{H(i)})).
\ee
Moreover, if $W_{0,i}-B_{H(i)}>W_{0,i'}-B_{H(i')}$ for all $i'\neq i$ and
$W_{0,i}-B_{H(i)}>0$, then every maximum weight matching
contains edge $(0,i)$. If $W_{0,i}-B_{H(i)}<0$ for all $i=1,\ldots,m$, then every maximum weight matching does not contain any
edge incident to $0$.
\end{lemma}

\begin{proof}
The proof is very similar to the one of Lemma \ref{lemma:recursionIndSet}.
If a maximum weight matching contains edge $(0,i)$ then its weight is $W_{0,i}$ plus maximum weight
of a matching in $H(i)$ with node $i$ excluded, plus the sum of the maximum weights of matchings in $T_{H(j)}, 1\leq j\leq m, j\neq i$.
If a maximum weight matching  contains none of the edges $(0,i), 1\leq i\leq m$, then its weight is simply sum of maximum
weight of matchings in $H(i)$ for all $i=1,\ldots,m$. We obtain
\be{eq:recursion2}
{\cal M}_H=\max(\max_{1\leq i\leq m}(W_{0,i}+{\cal M}_{H(i)}-B_{H(i)}+\sum_{j\neq i}{\cal M}_{H(j)}),\sum_{1\leq i\leq m}{\cal M}_{H(j)}).
\ee
Subtracting $\sum_{1\leq j\leq m}{\cal M}_{H(j)}$ from both sides of (\ref{eq:recursion2}) we obtain (\ref{eq:recursionM}). The proof
of the second part follows immediately from the discussion above. \qed
\end{proof}

\subsection{Long-range independence}\label{subsection:LongRangeIndependence}
We now consider trees $H$ of specific types. Given  integers $r\geq 3, d\geq 2$ let $H_r(d)$ denote an  $r$-regular finite
tree with depth $d$. The root node $0$ has degree $r-1$, all the nodes at  distance $\geq 1, \leq d-1$ from the root have
outdegree $r-1$, and all the nodes at distance $d$ from $0$ are leaves. (Usually, in the definition of
an $r$-regular tree, the root node is assumed to have degree $r$, not $r-1$.
The slight distinction here is done for convenience.)
Also, given a constant $c>0$,  a Poisson tree $H(c,d)$ with
parameter $c$ and depth $d$ is constructed as follows. The root node has a degree  which is a random
variable distributed according to $\Pois(c)$ distribution. All the children of $0$ have outdegrees which are also
random, distributed according to $\Pois(c)$. In particular, the children of $0$ have total degrees $1+\Pois(c)$.
Similarly, children of children of $0$ also have outdegree $\Pois(c)$, etc. We continue this process until
either the process stops at some depth $d'<d$, where no nodes in level $d'$ have any children, or until we reach
level $d$. In this case all the children of the nodes in level $d$ are deleted and the nodes in level $d$ become leaves. We obtain
a tree with depth $\leq d$. We call this a depth-$d$ Poisson tree.

Let $H=H_r(d)$ or $H(c,d)$. Suppose  the nodes and the edges of $H$ are equipped with weights $W_i, W_{i,j}$, which
are generated at random independently using a distribution function $F_w$.
Fix any infinite sequences $\bar w=(w_1,w_2,\ldots)\in [0,\infty)^{\infty}$
and $\bar b=(b_1,b_2,\ldots)\in \{0,1\}^{\infty}$. For every $i=1,2,\ldots,d$ let $i1,i2,\ldots,ij_i$ denote
the nodes of $H$ in level $i$ (if any exist for $H(c,d)$). When $H=H_r(d), j_i=(r-1)^i$, of course.
Let $({\cal I}|(\bar b,\bar w))$ denote the maximum weight of an independent set $V$
in $H$ such that the nodes $dj$ with $b_j=1$ are conditioned  to be in $V$, nodes $dj$ with $b_j=0$ are conditioned not
to be in $V$, and the weights of nodes $dj$ are conditioned to be equal to $w_j$ for  $j=1,\ldots,j_d$.
That is we are looking for maximum weight of an independent set among those which contain depth $d$
leaves with $b_j=1$, do not contain depth $d$ leaves with $b_j=0$, and with the weights of the leaves
deterministically set by $\bar w$. For brevity we call it  the maximum weight of an independent set with
boundary condition $(\bar b,\bar w)$. For the case  $H=H(c,d)$, the boundary condition is simply absent
when the tree does not contain any nodes in the last level $d$.
$({\cal I}_{H(ij)}|(\bar b,\bar w))$ are defined similarly for the subtrees $H(ij)$ spanned by nodes $ij$ in level $i$: given again
$\bar b,\bar w$, let $({\cal M}|(\bar b,\bar w))$ and  $({\cal M}_{H(ij)}|(\bar b,\bar w))$ denote, respectively, the maximum weight of a matching $E$
in $H$ and $H_{ij}$, such that the edges incident to nodes $dj$ are conditioned  to be in $E$ when $b_j=1$, edges incident to nodes $dj$
are conditioned not to be in $D$ when with $b_j=0$, and the weights of the edges incident to
nodes $dj$ are conditioned to be equal to $w_j, j=1,\ldots,j_d$ (of course, we refer to edges between nodes
in levels $d-1$ and $d$ as  there is only one edge per each node in level $d$).

For the case of independent sets, let $(B|(\bar b,\bar w))$ denote the bonus of the root node $0$ given the boundary condition $(\bar b,\bar w)$. Namely,
\[
(B|(\bar b,\bar w))=({\cal I}|(\bar b,\bar w))-\sum_{1\leq j\leq j_1}({\cal I}_{H(1j)}|(\bar b,\bar w)).
\]
For the case of matchings, let $(B|(\bar b,\bar w))$ also denote the bonus of the root node $0$ given the boundary condition $(\bar b,\bar w)$. Namely,
\[
(B|(\bar b,\bar w))=({\cal M}|(\bar b,\bar w))-\sum_{1\leq j\leq j_1}({\cal M}_{H(1j)}|(\bar b,\bar w)).
\]
It should be always clear from the context whether $B$ is taken with respect to independent sets or matchings.

The following theorem establishes the crucial \emph{long-range independence} property for the maximum weight independent
sets and matchings in trees $H=H_r(d),H(c,d)$ when the corresponding operator $T^2$ has a unique fixed point.
It establishes that the distribution of the bonus  $(B|(\bar b,\bar w))$ of the root
asymptotically is independent of the boundary condition $(\bar b,\bar w)$ as $d$ becomes large. Recall
our convention that $0$ denotes the distribution of a random variable which is zero to one with probability one.

\begin{theorem}\label{theorem:long range independence}
Given a  distribution function $F_w$ and  a regular tree $H=H_r(d)$ let $T=T_{{\cal I},r}$. Then for every $t\geq 0$ and $\bar b,\bar w$
\be{eq:LongRangeIndIndrOdd}
T^{d-1}(0)(t)\leq \pr((B|(\bar b,\bar w))\leq t)\leq T^{d-1}(F_w)(t)
\ee
when $d$ is odd, and
\be{eq:LongRangeIndIndrEven}
T^{d-1}(F_w)(t)\leq \pr((B|(\bar b,\bar w))\leq t)\leq T^{d-1}(0)(t)
\ee
when $d$ is even. Suppose in addition the operator $T^2=T^2_{{\cal I},r-1}$ has the unique fixed distribution $F^*$. Then
\be{eq:LongRangeIndIndr}
\sup_{\bar b,\bar w}\Big|\pr((B|(\bar b,\bar w))\leq t)-F^*(t)\Big|\rightarrow 0,
\ee
as $d\rightarrow\infty$. Similar assertion holds for $T=T_{{\cal M},r}$ and for $H=H(c,d)$ with
$T=T_{{\cal I},c}$ and $T=T_{{\cal M},c}$. For the cases $T=T_{{\cal M},r}$ and $T_{{\cal M},c}$, $F_w$
is replaced with $F_{w,r}$ and $F_{w,c}$ respectively.
\end{theorem}

Before we prove the proposition above, which is essentially the key result of the paper, let us compare it
with the developments in (\cite{AldousSteele:survey}). In that paper maximum weight matching is considered
on an  $n$-node tree, drawn independently and uniformly from the space of all  $n^{n-2}$ labelled trees. The notion of
a bonus is introduced and the recursion (\ref{eq:recursionM}) is derived. However, since a tree structure is assumed to begin
with, there is no need to consider the boundary conditions $(\bar b,\bar w)$. Here we avoid the difficulty of
the non-tree structure by proving the long-range independence property via the uniqueness of  fixed points of $T^2$.

\begin{proof}
We prove (\ref{eq:LongRangeIndIndrOdd}) and (\ref{eq:LongRangeIndIndrEven}) for independent sets in $H=H_r(d)$.
The proof of (\ref{eq:LongRangeIndIndr}) will then be almost immediate.
The proofs for other cases is similar. We will indicate the differences
where appropriate. The proof of (\ref{eq:LongRangeIndIndrOdd}) proceeds by induction in $d$. Suppose $d=1$. Then we have
\[
(B|(\bar b,\bar w))=\max(0,W_0-\sum_{i:1\leq i\leq r-1, \bar b_{1i}=1}\bar w_{1i})\leq W_0.
\]
Then $(B|(\bar b,\bar w))\prec F_w$. Trivially $0\prec (B|(\bar b,\bar w))$. This establishes the
bound for the case $d=1$. Suppose the assumption holds for $d-1$ and $d$ is odd. We have
\be{eq:B|bw}
(B|(\bar b,\bar w))=\max(0,W_0-\sum_{i:1\leq i\leq r-1}(B_{1i}|(\bar b,\bar w))).
\ee
By the inductive assumption the distribution of each of $(B_{1i}|(\bar b,\bar w))$ dominates
$T^{d-2}(F_w)$ and is dominated by $T^{d-2}(0)$. Applying this to (\ref{eq:B|bw}) we obtain
(\ref{eq:LongRangeIndIndrOdd}). The case $d$ is even is considered similarly. This completes the
induction and establishes (\ref{eq:LongRangeIndIndrOdd}) and (\ref{eq:LongRangeIndIndrEven}).
Now if $T^2$ has the unique fixed point distribution $F^*$ then by Proposition \ref{prop:Tconvergence}
$T^d(0),T^d(F_w)$ converge to $F^*$ as $d\rightarrow\infty$, and  we obtain (\ref{eq:LongRangeIndIndr}).
The proof for the case of Poisson tree $H(c,d)$ instead of
$H_r(d)$ is pretty much identical. For the proof of the same result for maximum weight matching we just use $F_{w,r}$ and $F_{w,c}$
instead of $F_w$, just as we did in the proof of Proposition \ref{prop:Tconvergence}. \qed
\end{proof}

While it is not important for the further results in this paper, it is interesting that the uniqueness of the
solution $T^2(F)=F$ is the tight condition for (\ref{eq:LongRangeIndIndr}), as the following theorem indicates.

\begin{theorem}\label{theorem:long range dependence}
Suppose the operator $T^2$ has more than one  fixed point distributions $F^*$. Then for every such $F^*$
\be{eq:LongRangeDepIndr}
\liminf_d\sup_{\bar b,\bar w}\Big|\pr((B|(\bar b,\bar w))\leq t)-F^*(t)\Big|>0.
\ee
\end{theorem}

\begin{proof}
As usual we start with $T=T_{{\cal I},r}$. The proofs for other cases are similar, we highlight the differences
where appropriate. Let $F_{**}$ and $F^{**}$ be distributions introduced in Proposition \ref{prop:Tconvergence}.
The non-uniqueness of the fixed point of $T^2$ implies using Proposition \ref{prop:Tconvergence} that $F_{**}\neq F^{**}$.
For every node $j$ in layer $d$ (last layer) of $H_r(d)$ set $b_j=0, w_j=0$. In particular, the bonus $B_j$
of each such node is zero. Then for every node in layer $d-1$ its bonus is given by the distribution $T(0)=F_w$,
the bonus of each node in layer $d-2$ has distribution $T^2(0)$, etc. The root node $0$ has bonus with
distribution $T^d(0)$. When $d$ is an even number diverging to infinity, from Proposition \ref{prop:Tconvergence},
$T^d(0)$ converges weakly to $F_{**}$. Thus the distribution of $(B|(\bar b,\bar w))$ converges to $F_{**}$.
If $F^*$, a fixed point of $T^2$, is distinct from $F_{**}$, then we obtain that  (\ref{eq:LongRangeDepIndr}) holds.
Suppose, on the other hand, $F^*=F_{**}$. We claim that $T(F_{**})\neq F_{**}$. Assuming this is the case
we consider the same boundary condition $\bar b=\bar w=0$ but take $d$ to be odd integer diverging to infinity.
Then the bonus of the root $0$ converges in distribution to $T(F_{**})\neq F_{**}$ and (\ref{eq:LongRangeDepIndr})
is shown again.

Assume $T(F_{**})=F_{**}$. Recall from the first part of Proposition \ref{prop:Tconvergence} that for every distribution
$F_0$, $T^{2s+1}(F_0)\prec T^{2s+1}(0)$. Taking $F_0=T(F^{**})$ we obtain $F^{**}=T^{2s+2}(F^{**})\prec T^{2s+1}(0)$.
Taking $s\rightarrow\infty$, $F^{**}\prec T(F_{**})=F_{**}$.  But from Proposition \ref{prop:Tconvergence} $F_{**}\prec F^{**}$,
and, as a result $F_{**}=F^{**}$ implying (again using Proposition \ref{prop:Tconvergence}) $T^2$ has the unique fixed point.
We obtained a contradiction. \qed
\end{proof}

\section{Applications to maximum weight independent sets and matchings in $G_r(n)$ and $G(n,c/n)$}\label{section:applications to G}
\subsection{Long-range independence in $G_r(n), G(n,c/n)$}\label{subsection:lop}
The goal of the current section is to demonstrate that Theorem \ref{theorem:long range independence} allows us to reduce the computation
of the maximum weight  independent set and the maximum weight matching in random graphs to a much simpler problem of finding those
in trees. We highlight this key message of the paper as the following
\emph{local optimality} property: if the operator $T^2$ corresponding to a maximum weight combinatorial object (independent set or matching)
in a sparse random graph has a unique fixed point, then for a randomly selected node (edge) of the graph, the event "the node (edge)
belongs to the optimal object" and the distribution of the node (edge) weight, conditioned that it does, asymptotically depends only on the
constant size neighborhood of the node and is independent from the rest of the graph. In other words, when $T^2$ has a unique fixed point,
the maximum weight independent sets and matchings exhibit a long-range independence property.

Our hope is that similar local optimality  can be established for other random combinatorial structures.

\begin{proof} [Proof of Theorem \ref{theorem:main}].
Again we start by proving the result for ${\cal I}_w(n,r)$. The proofs for other three objects is similar and we highlight
some differences in the end. Let $V_r\subset G_r(n)$ denote the  independent set which achieves the maximum weight ${\cal I}_w(n,r)$.
It is unique by continuity of the distribution $F_w$.
Consider a randomly selected node of the graph $G$, which, w.l.g., we may assume is node $0$.
By symmetry we have $\E[{\cal I}_w(n,r)]=n\E[W_01\{0\in V_r\}]$.
Let us fix a large positive integer $d$, which is a constant independent from $n$, and let $H(d)$ denote the depth-$d$
neighborhood of $0$. That is $H$ is the collection of nodes in $G$ which are connected to $0$ by paths with length $\leq d$.
It is well known that, w.h.p. as $n\rightarrow\infty$, $H$ is a depth-$d$ $r$-regular tree, \cite{JansonBook},
except in this case, unlike in Subsection \ref{subsection:LongRangeIndependence}, the root node has outdegree $r$
and the remaining non-leaf nodes have outdegree $r-1$.
Let $\partial H$ denote the leaves of this tree (level $d$). Fix any binary vector $\bar b$ with dimension
$|\partial H|+|G\setminus H|$, and any non-negative vector $\bar w$   also with dimension
$|\partial H|+|G\setminus H|$. Assume the vector $\bar b$ is such that if two nodes in $\partial H\cup (G\setminus H)$
are connected by an edge, only one of these two nodes can have the corresponding component of  $\bar b$ equal to $1$.
That is, the nodes marked $1$ by $\bar b$ correspond to some independent set in $\partial H\cup(G\setminus H)$.
Consider the problem of finding the maximum weight independent set $V_r$ in $G$ when the weights of the nodes in
$\partial H\cup (G\setminus H)$ are conditioned to be $\bar w$, nodes $i\in \partial H\cup (G\setminus H)$ with
the corresponding component of $\bar b$ equal to one are conditioned to belong to $V_r$ and the remaining
nodes in $\partial H\cup (G\setminus H)$ are conditioned not to belong to $V_r$. Notation-wise, we
consider the value $({\cal I}_w(n,r)|(\bar b,\bar w))$. In particular, we need to select to maximum
weight independent set in the tree $H$, which is consistent with conditioning $(\bar b,\bar w)$. Naturally,
the consistency needs to be checked only across the boundary $\partial H$.

Let $(B|(\bar b,\bar w))$ and $(B_i|(\bar b,\bar w)), 1\leq i\leq r$ denote the bonus of the node $0$ and
the bonuses of its neighboring nodes $1,2,\ldots,r$, respectively. By Lemma \ref{lemma:recursionIndSet}, we
have $(B|(\bar b,\bar w))=\max(0,W_0-\sum_{1\leq i\leq r}(B_i|(\bar b,\bar w)))$.
Since the distribution function $F_w$ of
$W_0$ is continuous, then $F^*$ is continuous as well and, as a  result, $W_0-\sum_{1\leq i\leq r}(B_i|(\bar,\bar w))=0$ with probability
zero. Then applying the second part of Lemma \ref{lemma:recursionIndSet}, node $0$ belongs to the maximum weight independent set if and only if
$B=W_0-\sum_{1\leq i\leq r}(B_i|(\bar b,\bar w))>0$. Therefore
\[
{\E[{\cal I}_w(n,r)]\over n}=\sum_{\bar b,\bar w}\E[W_01\{W_0-\sum_{1\leq i\leq r}(B_i|(\bar b,\bar w))>0\}]\pr(\bar b,\bar w).
\]
But by Theorem \ref{theorem:long range dependence}, for every $\bar b,\bar w$ the distribution of $(B_i|(\bar b,\bar w))$
converges to $F^*$, the unique fixed point of $T^2_{{\cal I},r-1}$ as $d$ becomes large. We conclude
\[
\lim_n{\E[{\cal I}_w(n,r)]\over n}=\E[W_01\{W_0-\sum_{1\leq i\leq r}B_i>0\}],
\]
where $B_i\sim F^*$.
This completes the proof for the maximum weight independent set in $G_r(n)$. When the graph $G(n,c/n)$ is considered the proof is very similar, we just
use the fact that $H$ -- depth $d$ neighborhood of $0$ approaches in distribution a Poisson tree \cite{TenLectures}. The proofs for maximum
weight matchings are similar, we use Lemma \ref{lemma:recursionM} instead of Lemma \ref{lemma:recursionIndSet}. \qed
\end{proof}

\begin{proof} [Proof of Corollary \ref{coro:main}]
In fact we have proved this result already en route of proving Theorem \ref{theorem:main} above. We have shown
that given  $\bar b,\bar w$ and conditioned on $H$ being a tree a fixed node $0$ belongs to the maximum
weight independent set iff $(B|(\bar b,\bar w))=W_0-\sum_{1\leq i\leq r}(B_i|(\bar b,\bar w))>0$. Repeating the proof of  Theorem \ref{theorem:main}
${\cal IN}_w(n,r)/n=\E[1\{0\in V_r\}]\rightarrow \E[1\{W_0-\sum_{1\leq i\leq r}B_i\}]$, where $B_i\sim F^*$. \qed
\end{proof}

\begin{proof}[Proof of Theorem \ref{theorem:Gaussian}]
The essential ingredients for the proof of this results are already established in the proof of Theorem \ref{theorem:main}
above. We start with ${\cal I}_w(n,r)$. The proofs for other cases are very similar and we omit them.
Let ${\cal O}={\cal O}(n,r)\subset [n]$ be the independent set achieving the
maximum weight. Fix two nodes $i,j\in [n]$ and arbitrary $\epsilon>0$.
Let  $H=H(i,d)$ denote the depth-$d$ graph-theoretic neighborhood of
$i$ in the graph $G=G_r(n)$ and let $\partial H$ denote the boundary of $H$ -- the nodes of $H$ at distance $d$ from $i$.
Let ${\cal E}_T$ denote the event that $H$ is a ($r$-regular) tree. From the theory of random regular graphs \cite{JansonBook},
$\pr({\cal E}_T)\rightarrow 1$ as $n\rightarrow\infty$. As a result $\pr({\cal E}_T)\geq 1-\epsilon$ for all $n\geq n_0$
for some $n_0=n_0(d)$ (note the dependence on $d$).
Fix any realization of $(G\setminus H)\cup \partial H$ together with the realization
of the weights $\bar w$ of nodes in $G\setminus H\cup \partial H$ and indicators $\bar b$ of whether the nodes belong
to ${\cal O}$. As far as deciding which  nodes of $H\setminus \partial H$ are in ${\cal O}$ and in particular whether node $i$ belongs to
${\cal O}$ only the restriction of $\bar w,\bar b$ to $\partial H$ is relevant. Conditioning on the event ${\cal E}_T$
denote by $(B|(\bar b,\bar w))$ the bonus of $i$ in $H$ (for completeness define $(B|(\bar b,\bar w))$ to be zero when
the event ${\cal E}_T$ does not hold).
Applying Lemma \ref{lemma:recursionIndSet} and the continuity of $F_w$,
we have $i\in {\cal O}$ iff $(B|(\bar b,\bar w))>0$. Applying Theorem \ref{theorem:long range independence}
$|\pr((B|(\bar b,\bar w))>0)-\pr(B>0)|<\epsilon$ for all $d\geq d_0(\epsilon)$ for some $d_0(\epsilon)$, and for any $\bar w,\bar b$, where $B\sim F^*$
and $F^*$ is the unique fixed point of $T^2=T^2_{{\cal I},r-1}$. Thus
\begin{eqnarray}
|\pr(i\in {\cal O})-\pr(B>0)| &\leq&  \sum_{\hat G\bar w,\bar b}|\pr(i\in {\cal O}|{\cal E}_T,\hat G,\bar w,\bar b)-\pr(B>0)|\pr(\hat G,\bar w,\bar b) \notag \\
&+& |\pr(i\in {\cal O}|\bar {\cal E}_T)-\pr(B>0)|\pr(\bar{\cal E}_T) \notag \\
&\leq& 2\epsilon, \label{eq:iinO}
\end{eqnarray}
whenever $d\geq d_0(\epsilon)$ and $n\geq n_0(d)$, where $\hat G$ denotes generically a realization of the subgraph $(G\setminus H)\cup \partial H$.
Observe that since $|H|\leq 1+r+(r-1)^2+\cdots+(r-1)^d$ then $\pr(j\notin H)\geq 1-\epsilon$ for all $n\geq n_1(d)$ for some $n_1(d)$.
Then
\begin{eqnarray}
&&|\pr(i,j\in{\cal O})-\pr(i\in{\cal O})\pr(j\in{\cal O})| \notag\\
&\leq & |\pr(i\in{\cal O}|{\cal E}_T,j\notin H,j\in{\cal O})\pr({\cal E}_T,j\notin H,j\in{\cal O})-\pr(i\in{\cal O})\pr(j\in{\cal O})| \label{eq:jnotinH}\\
&+& \pr(i\in{\cal O},{\cal E}_T,j\in H\cap {\cal O})+\pr(i\in{\cal O},\bar {\cal E}_T,j\in {\cal O}) \label{eq:jinH}
\end{eqnarray}
The event $j\notin H,j\in {\cal O}$ is completely described by the realizations $\hat G,\bar w,\bar b$. Since
$|\pr((B|(\bar b,\bar w))>0)-\pr(B>0)|<\epsilon$ for all $d\geq d_0(\epsilon)$ then
\be{eq:EjnotinH}
|\pr(i\in{\cal O}|{\cal E}_T,j\notin H,j\in{\cal O})-\pr(B>0)|<\epsilon,
\ee
for all $d\geq d_0(\epsilon)$. Also
$\pr(j\in {\cal O})=\pr({\cal E}_T,j\notin H,j\in{\cal O})+\pr({\cal E}_T,j\in H\cap {\cal O})+\pr(\bar {\cal E}_T,j\in{\cal O})$.
But $\pr({\cal E}_T,j\in H\cap {\cal O})\leq \pr(j\in H)\leq \epsilon$ for all $n\geq n_1(d)$
and $\pr(\bar {\cal E}_T,j\in{\cal O})\leq \pr(\bar {\cal E}_T)\leq \epsilon$ for all $n\geq n_0(d(\epsilon))$. As a result
\be{eq:O}
|\pr(j\in {\cal O})-\pr({\cal E}_T,j\notin H,j\in{\cal O})|\leq 2\epsilon
\ee
whenever $n\geq \max(n_0(d(\epsilon)),n_1(d(\epsilon)))$. Combining (\ref{eq:EjnotinH}), (\ref{eq:iinO}) and (\ref{eq:O}) we obtain
\begin{eqnarray*}
&&\pr(i\in{\cal O}|{\cal E}_T,j\notin H,j\in{\cal O})\pr({\cal E}_T,j\notin H,j\in{\cal O})-\pr(i\in{\cal O})\pr(j\in{\cal O}) \\
&\leq& (\pr(i\in{\cal O})+3\epsilon)(\pr(j\in{\cal O})+2\epsilon)-\pr(i\in{\cal O})\pr(j\in{\cal O}) \\
&\leq& 5\epsilon+\epsilon^2<6\epsilon
\end{eqnarray*}
Similarly we show
\[
\pr(i\in{\cal O}|{\cal E}_T,j\notin H,j\in{\cal O})\pr({\cal E}_T,j\notin H,j\in{\cal O})-\pr(i\in{\cal O})\pr(j\in{\cal O})\geq -6\epsilon.
\]
We conclude that the value in (\ref{eq:jnotinH}) is bounded by $6\epsilon$. Each summand in (\ref{eq:jinH}) is bounded by $\epsilon$
since $\pr(j\in H)$ when $n\geq n_1(d(\epsilon))$ and $\pr(\bar {\cal E}_T)\leq \epsilon$ when $n\geq n_0(d(\epsilon))$. We conclude
that whenever $n\geq \max(n_0(d(\epsilon)),n_1(d(\epsilon)))$, $|\pr(i,j\in{\cal O})-\pr(i\in{\cal O})\pr(j\in{\cal O})|\leq 8\epsilon$.
This completes the proof of the theorem. \qed
\end{proof}

\subsection{Computation of limits. Exponentially distributed weights}\label{subsection:ComputationLimitsExp}
We prove Theorem \ref{theorem:mainEXPr} in this subsection. Thanks to Theorem \ref{theorem:main}, we can
focus on proving the uniqueness and computing the fixed points of the operator $T^2$. As usual we start
with maximum weight independent set in $G_r(n)$. The analysis of other cases is similar and will follow
immediately. The calculations are similar to the ones in \cite{AldousSteele:survey} performed for maximum weight matching
in random trees. The  difference is that we have to compute the fixed point of the operator $T^2$ and not just $T$.

\vspace{.1in}

\begin{proof}[Proof of Theorem \ref{theorem:mainEXPr}]

\begin{itemize}
\item {\bf Independent sets in} $G_r(n)$.
Let $F^*$ denote any fixed point distribution of $T^2=T^2_{{\cal I},r-1}$ (at least one exists by Proposition \ref{prop:Tconvergence}),
and let $B\sim F^*$. Then $B=\max(W-\sum_{1\leq i\leq r-1}\hat B_i), W\sim \Exp(1), \hat B_i\sim T(F^*)$. Similarly,
if $\hat B\sim T(F^*)$, then $\hat B=\max(W-\sum_{1\leq i\leq r-1}B_i), W\sim \Exp(1), B_i\sim F^*$.
Then for any $t\geq 0$
\[
\pr(B>t)=\pr(W>\sum_{1\leq i\leq r-1}\hat B_i+t)=e^{-t}\pr(W>\sum_{1\leq i\leq r-1}\hat B_i),
\]
and similarly
\[
\pr(\hat B>t)=e^{-t}\pr(W>\sum_{1\leq i\leq r-1}B_i),
\]
where we use the memoryless property of the exponential distribution. Let $b=\pr(B=0), \hat b=\pr(\hat B=0)$
where $B\sim F^*$ and $\hat B\sim T(F^*)$. Our next goal is computing $b$ and $\hat b$. From above we obtain
\be{eq:B>t}
\pr(B>t)=e^{-t}(1-b), \qquad \pr(\hat B>t)=e^{-t}(1-\hat b),
\ee
implying
\be{eq:B>tConditioned}
\pr(B>t|B>0)\sim \Exp(1), \qquad \pr(\hat B>t|\hat B>0)\sim \Exp(1).
\ee
Then
\be{eq:b}
b=\pr(W-\sum_{1\leq i\leq r-1}\hat B_i\leq 0)=\int_0^{\infty}e^{-t}\pr(\sum_{1\leq i\leq r-1}\hat B_i\geq t)dt
\ee
In order to compute $\pr(\sum_{1\leq i\leq r-1}\hat B_i\geq t)$ we condition on $j\leq r-1$ terms $\hat B_i$
out of $r-1$ being equal to zero, and the rest positive. This occurs with probability $\binom{r-1}{j}\hat b^j(1-\hat b)^{r-1-j}$.
When $j<r-1$, the sum of $r-1-j$ non-zero terms $\hat B_i$ has an  Erlang distribution with parameter $r-1-j$ (sum of $r-1-j$ independent random
variables distributed as $\Exp(1)$). The density function of this distribution is $f(z)={z^{r-2-j}\over (r-2-j)!}e^{-z}$ and
the tail probability is
\be{eq:erlang}
\pr(~\cdot~>t)=\int_t^{\infty}{z^{r-2-j}\over (r-2-j)!}e^{-z}dz=\sum_{0\leq i\leq r-2-j}{t^i\over i!}e^{-t}
\ee
and
\[
\int_0^{\infty}e^{-t}\pr(~\cdot~>t)dt=\sum_{0\leq i\leq r-2-j}{1\over 2^{i+1}}=1-{1\over 2^{r-1-j}}.
\]
Combining with (\ref{eq:b}) and interchanging integration and summation, we obtain
\be{eq:brecursion}
b=\sum_{0\leq j\leq r-1}\binom{r-1}{j}\hat b^j(1-\hat b)^{r-1-j}(1-{1\over 2^{r-1-j}})=1-({1+\hat b \over 2})^{r-1}.
\ee
Similar calculations lead to $\hat b=1-({1+\hat b \over 2})^{r-1}$. Combining
\be{eq:ff}
b=f(f(b)), \qquad {\rm where} ~f(x)=1-({1+x\over 2})^{r-1}
\ee

\begin{lemma}\label{lemma:bhatb}
The equation (\ref{eq:ff}) has a unique solution $b^*$ within the range $b\in [0,1]$ iff $r\leq 4$.
In this case $b^*$ is also the unique solution of $f(b)=b$.
\end{lemma}

Figures 1 and 2 below show the graphs of $f(f(b))$
for the cases $r=4$ and $r=8$. The first corresponds to the case of the unique solution. In
the second case there are more than one solution.

\begin{proof}
Note that for every $r\geq 1$ the equation $x=f(x)$ has exactly one solution in $[0,1]$
since $f$ is a strictly decreasing function and $f(0)=1-1/2^{r-1}>f(1)=0$. This solution is also a solution to $x=f(f(x))$.

We now prove the uniqueness for $2\leq r\leq 4$ and non-uniqueness for $r>4$.
Let $r\leq 4$. We claim that for all $x\in [0,1]$.
\be{eq:ffx}
{df(f(x))\over dx}={(r-1)^2\over 4}\Big({1+f(x)\over 2}\Big)^{r-2}\Big({1+x\over 2}\Big)^{r-2}<1.
\ee
This would imply that $f(f(x))-x$ is a strictly decreasing function and is equal to
zero in at most one point. (\ref{eq:ffx}) is equivalent to
\be{eq:derivative<1}
(1+f(x))(1+x)<{2^{2r-2\over r-2}\over (r-1)^{2\over r-2}}.
\ee
Let us check the validity of this inequality at the end points $x=0,1$. Since $f(0)<1, f(1)=0$,
the left hand side in both points is at most $2$ which is strictly smaller than the  right-hand side
for $r=2,3,4$, as it is easily checked. It remains to check the inequality at the points where the
derivative of the function $g(x)=(1+f(x))(1+x)$ vanishes. We have
\begin{eqnarray*}
\dot g(x) & = & 1+\dot f(x)+f(x)+x\dot f(x) \\
& = & 1-{r-1\over 2}\Big({1+x\over 2}\Big)^{r-2}+1-\Big({1+x\over 2}\Big)^{r-1}-x{r-1\over 2}\Big({1+x\over 2}\Big)^{r-2} \\
& = & 2-\Big({1+x\over 2}\Big)^{r-1}-{(1+x)(r-1)\over 2}\Big({1+x\over 2}\Big)^{r-2} \\
& = & 2-r\Big({1+x\over 2}\Big)^{r-1}
\end{eqnarray*}
Thus $\dot g(x)=0$ in exactly one point $x=2({2\over r})^{1\over r-1}-1$. We check that (\ref{eq:derivative<1})
holds in this point for $r=2,3,4$.

We now prove the non-uniqueness of the solution to (\ref{eq:ff}) when $r\geq 5$. Let $b^*$
denote the unique fixed point of $f(b^*)=b^*$. Clearly $f(f(b^*))=b^*$. We claim that
\be{eq:b^*}
{df(f(x))\over dx}\Big|_{x=b^*}>1.
\ee
This implies the result since we get that for $\epsilon$ sufficiently small, $f(f(b^*-\epsilon))<b^*-\epsilon$. But
$f(f(0))>0$. Therefore there exists a different fixed point of (\ref{eq:ff}) in the interval $(0,b^*)$. To show (\ref{eq:b^*})
note that
\[
{df(f(x))\over dx}\Big|_{x=b^*}={(r-1)^2\over 4}\Big({1+f(b^*)\over 2}\Big)^{r-2}\Big({1+b^*\over 2}\Big)^{r-2}=
{(r-1)^2\over 4}\Big({1+b^*\over 2}\Big)^{2r-4}.
\]
We need to show that
\be{eq:br}
{(r-1)^2\over 4}\Big({1+b^*\over 2}\Big)^{2r-4}>1,
\ee
which is equivalent to $b^*>2({2\over r-1})^{1\over r-2}-1\equiv b(r)$. We claim that $b(r)<f(b(r))=1-({1+b(r)\over 2})^{r-1}$.
Since $b^*=f(b^*)$ and $f$ is a decreasing function, this would imply $b(r)<b^*$ or (\ref{eq:br}). Note, that $(1+b(r))/2=({2\over r-1})^{1\over r-2}$.
Thus we need to check that
\begin{eqnarray*}
2({2\over r-1})^{1\over r-2}-1<1-({2\over r-1})^{r-1\over r-2} & <=> & 2({2\over r-1})^{1\over r-2}<2-({2\over r-1})^{r-1\over r-2}  \\
& <=> & ({2\over r-1})^{1\over r-2}<1-{2^{1\over r-2}\over (r-1)^{1+{1\over r-2}}} \\
& <=> & ({2\over r-1})^{1\over r-2}<{r-1\over r}  \\
& <=> & 2<(r-1)(1-{1\over r})^{r-2}.
\end{eqnarray*}
Note that $g(r)=(1-{1\over r})^{r-2}$ is a strictly growing function of $r$ since $\log(g(r))=(r-2)\log(1-{1\over r})$. Therefore,
it suffices to verify that for $r=5$ we have $2<4(1-{1\over 5})^3\approx 2.048$. \qed
\end{proof}

\begin{figure}
\begin{center}
\scalebox{.5}{
\includegraphics{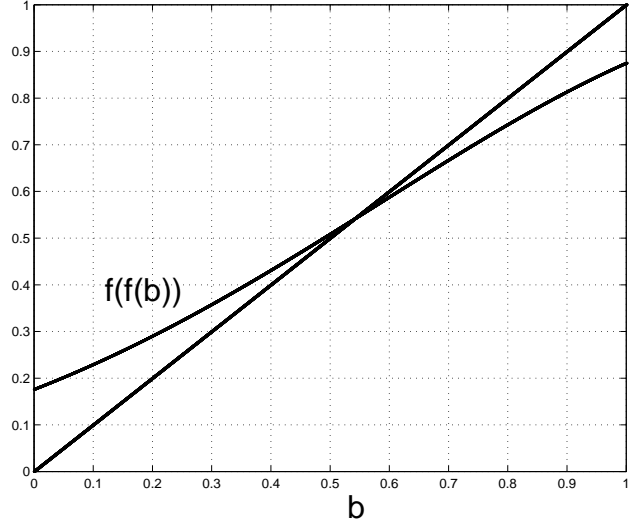}\label{r4}
}
\end{center}
\caption{$f(f(b))=b$ has one solution when $r=4$.}
\end{figure}

\begin{figure}
\begin{center}
\scalebox{.5}{
\includegraphics{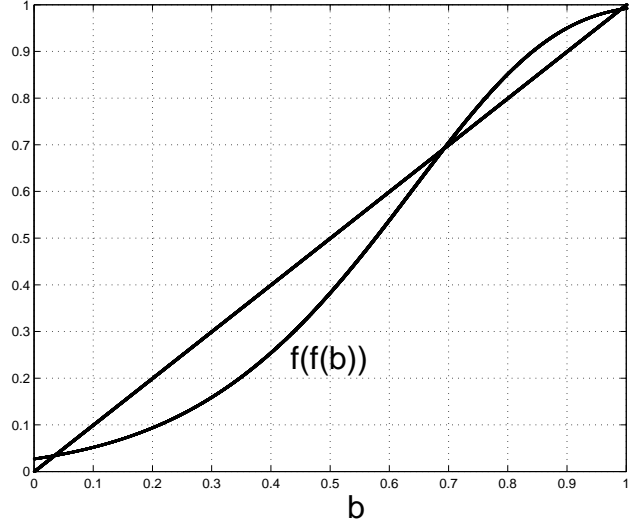}\label{r8}
}
\end{center}
\caption{$f(f(b))=b$ has more than one solution when $r=8$.}
\end{figure}

Applying Lemma \ref{lemma:bhatb} $b=\hat b$ is the unique solution of $f(b)=b$.
When $r=2$, we obtain $b=1/3$.
When $r=3$ we obtain $b=1-(1+b)^2/4$ or $b=2\sqrt{3}-3\approx .4641~$. When $r=4$ we find numerically that $b\approx .5419$.

Recall, that the distribution $F^*$ of $B$ is given by $\pr(B=0)=b$ and $\pr(B>t)=(1-b)e^{-t}$.
Applying (\ref{eq:mainIndr}) of Theorem \ref{theorem:main} and generating $B_1,\ldots,B_r\sim F^*$
independently, we obtain
\be{eq:lim0}
\lim_n {{\cal I}_w(n,r)\over n}=\int_0^{\infty}t e^{-t}\pr(t>\sum_{1\leq i\leq r}B_i)dt
\ee
We again condition on $j$ out of $r$ terms $B_i$ to be equal to zero and the rest positive.
This occurs with probability $\binom{r}{j}b^j(1-b)^{r-j}$. W.l.g. assume that the non-zero terms are $B_1,\ldots,B_{r-j}$.
We have
\be{eq:lim1}
\int_0^{\infty}t e^{-t}\pr(t>\sum_{1\leq i\leq r-j}B_i)dt=1-\int_0^{\infty}t e^{-t}\pr(\sum_{1\leq i\leq r-j}B_i>t)dt.
\ee
Repeating the calculations (\ref{eq:erlang}), we obtain that
$\pr(\sum_{1\leq i\leq r-j}B_i>t)=\sum_{0\leq i\leq r-1-j}{t^i\over i!}e^{-t}$. Then the expression of the right-hand side of
(\ref{eq:lim1}) becomes
\be{eq:lim2}
1-\int_0^{\infty}t e^{-t}\sum_{0\leq i\leq r-1-j}{t^i\over i!}e^{-t}dt=1-\sum_{0\leq i\leq r-1-j}{i+1\over 2^{i+2}}={r-j+2\over 2^{r-j+1}}.
\ee
Combining, we obtain
\[
\lim_n {{\cal I}_w(n,r)\over n}=\sum_{0\leq j\leq r}\binom{r}{j}b^j(1-b)^{r-j}{r-j+2\over 2^{r-j+1}}=
{1\over 2}2\sum_{0\leq j\leq r}\binom{r}{j}b^j({1\over 2}-{b\over 2})^{r-j}+
\]
\[
{1\over 2}\sum_{0\leq j\leq r}\binom{r}{j}b^j({1\over 2}-{b\over 2})^{r-j}(r-j)\Big).
\]
The first summand is simply $({1\over 2}+{b\over 2})^r$. We compute the second summand using the following
probabilistic argument. Rewrite the expression as
\[
({1\over 2}+{b\over 2})^r\sum_{0\leq j\leq r}\binom{r}{j}{b^j({1\over 2}-{b\over 2})^{r-j}\over ({1\over 2}+{b\over 2})^r}(r-j).
\]
The sum above is simply the number of successes in $r$  Bernoulli trials with the probability of success equal
to $({1\over 2}-{b\over 2})/({1\over 2}+{b\over 2})=(1-b)/(1+b)$. Namely, it is $r(1-b)/(1+b)$.
We obtain
\be{eq:Iwr}
\lim_n {{\cal I}_w(n,r)\over n}=({1\over 2}+{b\over 2})^r(1+{r(1-b)\over 2+2b})={(1+b)^{r-1}\over 2^{r+1}}(2+2b+r-r b).
\ee
Recall from Lemma \ref{lemma:bhatb} that $(1+b)^{r-1}/2^{r-1}=1-b$. This proves (\ref{eq:mainIndr=3})
Plugging the corresponding value of $b$ for $r=2,3,4$, we obtain $\lim_n {{\cal I}_w(n,2)\over n}={2\over 3}$,
$\lim_n {{\cal I}_w(n,3)\over n}=.6077$ and $\lim_n {{\cal I}_w(n,4)\over n}=.5632$. This concludes the proofs for the case of
maximum weight independent set in $G_3(n)$ and $G_4(n)$. Before we continue with other cases, it is convenient to prove
Corollary \ref{coro:mainCardinality}, as the proof is almost immediate from above.

\begin{proof}[Proof of Corollary \ref{coro:mainCardinality}]
We have established above that $T^2_{{\cal I},r}$ has a unique fixed point iff $r\leq 4$.
Applying (\ref{eq:mainIndrCard}) we need to compute $\E[1\{W-\sum_{1\leq i\leq r}B_i>0\}]=\pr(W>\sum_{1\leq i\leq r}B_i)$,
where $W\sim \Exp(1), B_i\sim F^*$. Instead of computing this quantity directly, note that the probability above
is exactly $1-b$, if we the summation above was up to $r-1$ not $r$. Repeating the computations up to (\ref{eq:brecursion}),
we obtain
\[
\pr(W>\sum_{1\leq i\leq r}B_i)=({1+b \over 2})^r.
\]

Plugging the obtained values of $b$ for $r=2,3,4$, we obtain (\ref{eq:mainIndrCard=2,3,4}). \qed
\end{proof}

\item {\bf Independent sets in} $G(n,c/n)$.

Let $F^*$ denote now any fixed point of $T^2=T^2_{{\cal I},c}$. We introduce again $b=F^*(0)=\pr(W-\sum_{1\leq i\leq m}\hat B_i\leq 0)$,
where $\hat B_i\sim T(F^*)$ and $m\sim \Pois(c)$. Similarly, $\hat b=T(F^*)(0)=\pr(W-\sum_{1\leq i\leq m}B_i\leq 0)$,
where $B_i\sim F^*$. Repeating the computations done for $G_r(n)$, we obtain similarly to (\ref{eq:brecursion}) that conditioning on $m=k$
\[
b=1-({1+\hat b \over 2})^k.
\]
Then $b=\sum_{k\geq 0}{c^k\over k!}e^{-c}(1-({1+\hat b \over 2})^k)=1-e^{-{c(1-\hat b)\over 2}}$. Similarly,
$\hat b=1-e^{-{c(1-b)\over 2}}$. Thus, $b$ must satisfy $1-b=\exp(-{c\over 2}\exp(-{c\over 2}(1-b)))$.
Recall, however, that by the second part of Theorem \ref{theorem:KarpSipser}, the equation above
has a unique solution iff $c\leq 2e$. In this case $b$ is also the unique solution of $1-b=\exp(-{c\over 2}(1-b))$.
We now apply (\ref{eq:mainIndr}) of Theorem \ref{theorem:main} to compute $\lim_n{{\cal I}_w(n,c)\over n}$, where we
substitute $m\sim \Pois(c)$ for $r$. In order to shortcut the computations, we use (\ref{eq:lim0}) and (\ref{eq:Iwr}).
We have
\[
\lim_n{{\cal I}_w(n,c)\over n}=\sum_{k\geq 0}{c^k\over k!}e^{-c}\int_0^{\infty}t e^{-t}\pr(t>\sum_{1\leq i\leq k}B_i)dt,
\]
where $B_i\sim F^*$. Recall, though, from (\ref{eq:Iwr}), that
\[
\int_0^{\infty}t e^{-t}\pr(t>\sum_{1\leq i\leq k}B_i)dt=({1+b\over 2})^k(1+{k(1-b)\over 2(1+b)}).
\]
Combining
\[
\lim_n{{\cal I}_w(n,c)\over n}=e^{-{c(1-b)\over 2}}+{(1-b)\over 2(1+b)}{c(1+b)\over 2}e^{-{c(1-b)\over 2}}=
(1+{c(1-b)\over 4})e^{-{c(1-b)\over 2}}.
\]
Using $b=1-\exp(-{c(1-b)\over 2})$, we obtain (\ref{eq:mainIndc}).

\item {\bf Matchings in } $G_r(n)$.

Let $F^*$ denote any fixed point distribution of $T=T_{{\cal M},r}$
Then we have the following distributional identities
$\hat B=\max_{1\leq i\leq r-1}(0,W_i-B_i), W_i\sim \Exp(1), B_i\sim F^*,\hat B\sim T(F^*)$ and
$B=\max_{1\leq i\leq r-1}(0,W_i-\hat B_i), W_i\sim \Exp(1), \hat B_i\sim T(F^*),B\sim F^*$.
Let  $\hat b=\pr(W-\hat B<0)$, where $W\sim \Exp(1)$ and $\hat B\sim T(F^*)$,
and let $b=\pr(W-B<0)$, where $B\sim F^*$. Then for any $t\geq 0$
\be{eq:BdistrM}
\pr(B\leq t)=\pr(\max_{1\leq i\leq r-1}(W_i-\hat B_i)\leq t)=(1-\pr(W_1>t+\hat B_1))^{r-1}=(1-e^{-t}(1-\hat b))^{r-1}.
\ee
In particular
\be{eq:bM}
\pr(B=0)=\hat b^{r-1}, \qquad d\pr(B\leq t)=(r-1)(1-\hat b)e^{-t}(1-e^{-t}(1-\hat b))^{r-2}dt, ~t>0
\ee

(We note as above that $\pr(W_1=\cdot)=0$, since $W_1$ has a continuous distribution). Then
\begin{eqnarray}
b & = & \int_0^{\infty}e^{-t}\pr(B>t)dt   \notag\\
& = &  \int_0^{\infty}e^{-t}(1-(1-e^{-t}(1-\hat b))^{r-1})dt \notag\\
& = & 1-\int_0^{\infty}e^{-t}(1-e^{-t}(1-\hat b))^{r-1}dt  \notag\\
& = & 1-\int_0^{\infty}e^{-t}(1-e^{-t}(1-\hat b))^{r-1}dt  \notag\\
& = & 1- \int_0^{\infty}(1-e^{-t}(1-\hat b))^{r-1}d(-e^{-t})  \notag\\
& = & 1- \int_0^1(1-z(1-\hat b))^{r-1}dz \notag\\
& = & 1- {(1-z(1-\hat b))^r\over r(-(1-\hat b))}\Big|_0^1 \notag \\
& = & 1-{1-\hat b^r\over r(1-\hat b)}. \label{eq:bM2}
\end{eqnarray}
Similarly, we obtain
\[
\hat b=1-{1-b^r\over r(1-b)},
\]
and combining we conclude that $b$ must be a fixed point of the
equation $f(f(x))=x$, where
\[
f(x)=1-{1-x^r\over r(1-x)}.
\]
\begin{lemma}\label{lemma:bhatbMatching}
For every $r\geq 2$ the equation $f(f(x))=x$ has a unique solution $x^*$ in the range $[0,1]$, which is the
unique solution of the equation $f(x)=x$.
\end{lemma}

\begin{proof}
Note that $(1-x^r)/(1-x)=1+x+\cdots+x^{r-1}$ and therefore $f(x)$ is a strictly decreasing function with $f(0)=1-1/r, f(0)=0$.
Therefore $f(x)=x$ has exactly one solution $x^*$. We now prove that $f(f(x))>x$ for all $x<x^*$ and $f(f(x))<x$ for all $x>x^*$.
This would complete the proof of the lemma. We need to show that for $x<x^*$
\[
f(f(x))=1-{1-f^r(x)\over r(1-f(x))}>x,
\]
which is equivalent to $1-f^r(x)<r(1-f(x))(1-x)$. But since $f$ is a decreasing function and $f(x^*)=x^*$, then
$f(x)>x$ for all $x<x^*$ and therefore $1-f^r(x)<1-x^r=r(1-f(x))(1-x)$.
Similarly, when $x>x^*$ we have $f(x)<x^*$ and then $1-f^r(x)>1-x^r=r(1-f(x))(1-x)$, resulting in
\[
f(f(x))=1-{1-f^r(x)\over r(1-f(x))}<x.
\]
\qed
\end{proof}

We conclude that $b=\hat b$ is determined as the unique solution of
\be{eq:bsolution}
b=1-{1-b^r\over r(1-b)},
\ee
and the unique fixed point of $T^2$ is the distribution given by (\ref{eq:BdistrM}) with
$b=\hat b$ given above.
Now, using (\ref{eq:mainMr}) of Theorem \ref{theorem:main} and (\ref{eq:bM})
we have

\begin{eqnarray}
\lim_n{{\cal M}_w(n,r)\over n} & = & \E[\sum_{1\leq i\leq r}W_i1\{W_i-B_i=\max_{1\leq j\leq r}(0,W_j-B_j)\}], \notag \\
& = & r\E[W_11\{W_1-B_1>\max_{2\leq j\leq r}(0,W_j-B_j)\}],  \notag \\
& = & r\int_0^{\infty}te^{-t}b^{r-1} \pr(t>\max_{2\leq j\leq r}(W_j-B_j))dt  \label{eq:B=0}\\
& + & r\int_0^{\infty}\int_0^tte^{-t}(r-1)(1-b)\times \notag \\
& \times & e^{-z}(1-e^{-z}(1-b))^{r-2}\pr(t-z>\max_{2\leq j\leq r}(W_j-B_j))dtdz,  \label{eq:B>0}
\end{eqnarray}
where the summands (\ref{eq:B=0}) and (\ref{eq:B>0}) corresponds to conditioning on $B_1=0$ and $B_1=z>0$. We now compute the integrals
in these summands. We have
\begin{eqnarray*}
\int_0^{\infty}te^{-t}b^{r-1} \pr(t>\max_{2\leq j\leq r}(W_j-B_j))dt & = & \int_0^{\infty}te^{-t}b^{r-1} \pr^{r-1}(t>W_2-B_2)dt \\
& = & \int_0^{\infty}te^{-t}b^{r-1} (1-\pr(W_2>t+B_2))^{r-1}dt \\
& = & \int_0^{\infty}te^{-t}b^{r-1} (1-e^{-t}\pr(W_2>B_2))^{r-1}dt \\
& = & \int_0^{\infty}te^{-t}b^{r-1} (1-e^{-t}(1-b)))^{r-1}dt  
\end{eqnarray*}
Similarly, we obtain that the integral  in (\ref{eq:B>0}) is equal to
\begin{eqnarray*}
& = & \int_0^{\infty}\int_0^tte^{-t}(r-1)(1-b)e^{-z}(1-e^{-z}(1-b))^{r-2}(1-e^{-t}(1-b)))^{r-1}dtdz \\
& = & \int_0^{\infty}\int_0^tte^{-t}(r-1)(1-b)(1-e^{-z}(1-b))^{r-2}(1-e^{-t}(1-b)))^{r-1}dtd(-e^{-z}) \\
& = & \int_0^{\infty}te^{-t}(r-1)(1-b)(1-e^{-t}(1-b)))^{r-1}dt\int_0^{e^{-t}}(1-w(1-b))^{r-2}dw \\
& = & \int_0^{\infty}te^{-t}(1-e^{-t}(1-b)))^{r-1}dt((1-w(1-b))^{r-1}\Big|^0_{e^{-t}}) \\
& = &  \int_0^{\infty}te^{-t}(1-e^{-t}(1-b)))^{r-1}dt- \int_0^{\infty}te^{-t}(1-e^{-t}(1-b)))^{2r-2}dt
\end{eqnarray*}
Substituting to summands in (\ref{eq:B=0}) and (\ref{eq:B>0}), we obtain (\ref{eq:mainMr=3}).

\item{\bf Matchings in } $G(n,c/n)$.

The derivation is very similar to the one for ${\cal M}_w(n,r)$. We introduce $b$ and $\hat b$
exactly as above. The equation (\ref{eq:bM2}) becomes
\begin{eqnarray*}
b & = & \sum_{m\geq 0}{c^m\over m!}e^{-c}(1-{1-\hat b^{m+1}\over (m+1)(1-\hat b)}) \\
  & = & 1-{1\over c(1-\hat b)}\sum_{m\geq 0}{c^{m+1}\over (m+1)!}e^{-c}+{1\over c(1-\hat b)}\sum_{m\geq 0}{c^{m+1}\hat b^{m+1}\over (m+1)!}e^{-c} \\
  & = & 1-{1\over c(1-\hat b)}+{e^{-c(1-\hat b)}\over c(1-\hat b)},
\end{eqnarray*}
Then $1-b$ is a fixed point of the equation $f(f(x))=x$, where $f(x)={1-e^{-cx}\over cx}$. That is
\[
{1-e^{c{1-e^{-cx}\over cx}-c}\over c{1-e^{-cx}\over cx}}=x.
\]
First we note that $x=0$ does not satisfy $f(f(x))=x$, so we assume $x>0$. The expression above then
becomes $e^{-cx}+cx^2-1=0$. The function $e^{-cx}+cx^2-1$ is a strictly convex function which is equal to zero
at $x=0$, has derivative $-c<0$ at $x=0$ and diverges to infinity as $x$ diverges to infinity.
Therefore its graph has exactly one  intersection with the horizontal axis in $x\in (0,\infty)$. Note
that at $x=1$ the value is $e^{-c}+c-1>0$ (for every $c>0$) therefore  there exists exactly solution to
$e^{-cx}+cx^2-1=0$ in $x\in (0,1)$. We conclude that  $b$ is uniquely
determined by the equation $e^{-c(1-b)}+cb^2-1=0$. The remainder of the calculation is done just like for $G_r(n)$
by conditioning first on specific values of $r$ and recalling $m=r$ with probability ${c^r\over r}!e^{-c}$. We omit
the fairly straightforward calculations. \qed
\end{itemize}
\end{proof}

\subsection{Extensions: Cycles and (non-random)  low degree graphs}\label{subsection:extensions}
In this section we first prove Corollary \ref{coro:cycle}. We then focus on proving Theorem \ref{theorem:WorstcaseGirth},
thereby improving the previous bound of Hopkins and Staton \cite{HopkinsStaton}.

\begin{proof}[Proof of Corollary \ref{coro:cycle}] The result is essentially established en route of proving Theorems \ref{theorem:main}
and \ref{theorem:mainEXPr}. Recall that in the proof of Theorem \ref{theorem:main} the only place we used the randomness of our
regular graph $G_r(n)$ was to say that a constant depth-$d$ neighborhood $H(d)$ of a randomly selected node $i\in[n]$ is a
depth-$d$  $r$-regular tree w.h.p. In the case of a cycle any constant depth-$d$ neighborhood of any node is w.p.1 a  path of length
$2d$ with the selected node in the middle. This is depth-$d$ $2$-regular tree. The answer for the maximum weight
matching on a cycle is the same as for independent set since we may simply assume the weights are assigned to nodes which
are to the left for each edge. \qed
\end{proof}

\begin{proof}[Proof of Theorem \ref{theorem:WorstcaseGirth}] As in the proof of Corollary \ref{coro:cycle}, we recall
that the only place the randomness of $G_r(n)$ was used in the proof of Theorem \ref{theorem:main} and Theorem \ref{theorem:mainEXPr}
was the fact that a constant depth-$d$ neighborhood of a randomly chosen node $i\in [n]$ is w.h.p. a depth $d$ $r$-regular tree.
In order to continue the proof we introduce a class of almost $r$-regular (non-random) graphs with large girth.
Given positive integers $d,R$ let ${\cal G}_r(n,d,R)$ denote the class of $n$-node graphs with girth at least $d$
and such that all but at most $R$ nodes have degree $r$ and the remaining nodes have degree less than $r$.
Given $G\in {\cal G}_r(n,d,R)$ a random node $i\in [n]$ its
depth-$d$ neighborhood $H(i,d)$ is w.h.p. (with respect to randomness of the choice of $i$) a depth-$d$ $r$-regular tree.
 Repeating the proofs of Theorem \ref{theorem:main} and \ref{theorem:mainEXPr}
and Corollary \ref{coro:mainCardinality} we conclude that for the case $r=3$ for every constant $R$ and for
all sufficiently large $n,d$, the expected maximum weight of an independent set in $G$
is given by (\ref{eq:mainIndr=2,3,4}) with $r=3$, that is $\approx.6077n$.
Also by Corollary \ref{coro:mainCardinality}, the expected cardinality  of the independent set achieving this weight is
given by (\ref{eq:mainIndrCard=2,3,4}) with $r=3$, that is $\approx .3923n$, also for all $n$ and $d$ sufficiently large.
Then the maximum cardinality of an independent set in $G$ is at least $.3923n-o(n)$  and this is a non-probabilistic bound
since we consider a maximum cardinality independent set in a fixed graph $G$. Thus for every constant  $R$
\[
\liminf_d\liminf_n\min_{G\in{\cal G}_3(n,d,R)}{{\cal I}(G)\over n}\geq .3923.. ~.
\]
To finish the proof we need to obtain similar lower bound for the class ${\cal G}(n,3,d)$ -- all graphs with degree
\emph{at most} $r=3$ and girth at least $d$. Suppose we are given any such $n$-node graph $G$.
Consider any two nodes $i_1,i_2$ which have degree $<3$ such that $i_2\notin H(i_1,d)$
(if any such  pair exist). Connect $i_1,i_2$ by an edge. This operation can only decrease the size of the largest
independent set in $G$. We claim that in addition the resulting graph still has girth at least $d$. In fact, if the
added edge becomes a part of a cycle with length less than $d$, then in the original graph the nodes $i_1,i_2$ were
connected by path of length at most $d-2$, contradicting the condition $i_2\notin H(i_1,d)$. Therefore the resulting
graph still has girth at least $d$. We continue this operation for every such pair $i_1,i_2$ until we either do not
have any nodes with degree less than $3$ at all or all such nodes belong to $H(i,d)$ for some node $i$. Since
$|H(i,d)|\leq 1+3+3^2+\ldots 3^d\equiv R$, we obtain a graph $G'\in{\cal G}_r(n,d,R)$. Then
\[
\liminf_d\liminf_n\min_{G\in{\cal G}(n,3,d)}{{\cal I}(G)\over n}\geq
\liminf_d\liminf_n\min_{G\in{\cal G}_3(n,d,R)}{{\cal I}(G)\over n}\geq .3923n.
\]
This completes the proof of the theorem. \qed
\end{proof}

\subsection{Deterministic and Bernoulli weights}\label{subsection:ComputationLimitsDet}
Are the results obtained above relevant to the case when the weight of each node
and edge $W=1$ with probability one or, in general, when the weights
take some discrete values? Let us examine these questions with respect to our usual four
operators $T$. We start with the case $W=1$.
For $T=T_{{\cal I},r}$, the corresponding distributional equation is $B=\max(0,1-\sum_{1\leq i\leq r-1}B_i)$.
If $B_i=0$, w.p.1, then $B=1$, w.p.1, and vice verse. Thus $B=0$ and $B=1$ are two fixed points of $T^2$, and
$T^2$ does not have a unique fixed point distribution. There is a physical explanation for the lack of uniqueness,
coming directly from the lack of long-range independence. Given a depth-$d$ $r$-regular tree $T$ with all the weights equal to unity,
note that the boundary does carry a non-vanishing information about the root in the following sense.
If all the leaves of the tree (boundary nodes) are conditioned to belong to the maximum weight independent set, then all the
parents of leaves cannot be part of the set. Then the maximum independent set is obtained by selecting all the nodes in level $d$,
not selecting nodes in level $d-1$, selecting all the nodes in level $d-2$,  and so on. In the end whether the root is selected is fully determined by the parity of $d$.
Thereby we do not have a long-range independence.
Contrast this with the discussion in Brightwell and Winkler \cite{BrightwellWinkler}, where
similar observation is used to show long-range dependence for Gibbs measures on infinite regular trees
for the hard-core model. It is not hard to see
that the similar lack of long-range independence holds for maximum weight matchings, where the weights are all $1$.

The situation is different for $T=T_{{\cal I},c}$. Let $F$ be a distribution function given
by $F(t)=p, t\in [0,1), F(1)=1$. Namely, $F$ is simply a Bernoulli distribution with parameter $p$ ($\Be(p)$).
If $B=\max(0,1-\sum_{1\leq i\leq m}B_i)$, where $B_i\sim F$ and $m\sim \Pois(c)$, then
$B=1$ if $\sum_{1\leq i\leq m}B_i=0$, which occurs with probability $\sum_{k\geq 0}{c^k\over k!}e^{-c}p^k=e^{-c(1-p)}$,
and $B=0$ otherwise. Thus $T(F)$ is $\Be(p_1)$ where $p_1=1-e^{-c(1-p)}$. Similarly $T^2(F)$ is $\Be(p_2)$, with
with  $p_2=1-e^{-ce^{-c(1-p)}}$. In general, for $s=1,2,\ldots$, $T^{2s}(F)$ is $\Be(p_{2s})$ with
$1-p_{2s}=e^{-ce^{-c(1-p_{2s-2})}}$. By Proposition \ref{prop:Tconvergence} we know that for $F=\Be(0)$ and $F=\Be(1)$,
$T^{2s}(F)$ converges
to some fixed point distributions $F_{**},F^{**}$  which by the argument above are $\Be(p_{**}),\Be(p^{**})$ with both
$p=p_{**}$ and $p=p^{**}$ satisfying
$1-p=e^{-ce^{-c(1-p)}}$. Recall from Theorem \ref{theorem:KarpSipser}, that the equation $x=e^{-ce^{-cx}}$
has a unique solution iff $c\leq e$. By Proposition \ref{prop:Tconvergence} this implies that when $c\leq e$, $T^s(F_0)$ converges
to $\Be(p^*)$ for any starting distribution $F_0$. Again applying Proposition \ref{prop:Tconvergence} we obtain  that $T^s(F_0)$ converges
to $\Be(p^*)$, $T$ and $T^2$ have the same unique fixed point distribution -- $\Be(p^*)$, where
$p^*$ is also the unique solution of $1-p^*=e^{-c(1-p^*)}$.

It is a simple exercise to see that the same holds for $T=T_{{\cal M},c}$: $T^2$ has  a unique fixed point iff $c\leq e$,
in which case the fixed point distribution is also $\Be(p^*)$. This is, of course, fully consistent
with Theorem \ref{theorem:KarpSipser}. We summarize these observations.

\begin{prop}\label{prop:fixedDeterministic}
Let the nodes and edges weights be equal to one, w.p.1. For every $r\geq 1$, $T^2_{{\cal I},r}, T^2_{{\cal M},r}$ have
at least two fixed point distributions. $T^2_{{\cal I},c}, T^2_{{\cal M},c}$ have the unique fixed point distribution
iff $c\leq e$ in which case  the unique fixed point distribution is $\Be(p^*)$, where $p^*$ is the unique solution of
$1-p^*=e^{-c(1-p^*)}$.
\end{prop}

Can we fully reproduce Theorem \ref{theorem:KarpSipser} for the case $c\leq e$? The problem with the case $F_w=1$
as, generally, with non-continuous distributions $F_w$, is that the probability of  $W-\sum_{1\leq i\leq r}B_i=0$
is no longer zero. As a result we do not have an exact condition purely in terms of $B=\max(0,W-\sum_{1\leq i\leq r}B_i)$
for whether the root node belongs to say maximum weight independent set. One natural approach would be to approximate
$F_w$ with a continuous distribution.
But the difficulty is the lack of closed form expression for the solution of fixed point of $T^2(F^*)=F^*$. Such a
solution $F^*$ can, though, also be approximated by computing $T^{2s}(0)$ and $T^{2s}(F_w)$ ($T^{2s}(F_{w,r}),T^{2s}(F_{w,c})$)
for matching) for large $s$ such that  the differences between the two distributions is sufficiently small.

Suppose now the weights $W_i$  of the nodes are  distributed as $\Be(z)$ for some parameter $z\in [0,1]$. For simplicity we will only consider
the case $r=3$ and $T=T_{{\cal I},3}$ and obtain a complete criteria for uniqueness.
\begin{prop}\label{prop:fixedBernoulli}
For $r=3$ the operator $T^2_{{\cal I},r-1}$ has a unique fixed point distribution iff $z\in [{1\over 4},1]$. In this case
the fixed point distribution is $\Be(p)$ with $p={\sqrt{5-4z}-1 \over 2(1-z)}$.
\end{prop}

\begin{proof} Let $F=\Be(p)$ for any $x\in [0,1]$. Then for $B'\sim T(F)$
we have $B'=\max(0,W-B_1-B_2)$ where $B_1,B_2$ are independent and distributed as $F$. Then
$B'=0$ when $W=0$ or when $W=1$ and $B_1+B_2>0$. This occurs with probability $z+(1-z)(1-p^2)=1-(1-z)p^2$.
Thus $B'\sim \Be(1-(1-z)p^2)$. Repeating the development above for the case of the deterministic weight,
we need to analyze the number of the solutions to the equation $f(f(x))=x$ where $f(x)=1-(1-z)x^2$. First,
the equation $f(x)=x$ leads to the unique solution $x^*={-1+\sqrt{5-4z} \over 2(1-z)}$ (it is simple to check that the solution
is in $[0,1]$ for every $z\in [0,1]$). Also $x^*$ is a solution of $g(x)\equiv f(f(x))=x$. We need to show that this is the unique solution
of this equation iff $z>1/4$.

First we show that when $z<1/4$, $\dot g(x)>1$. Since $g(0)>0$ this would imply that there exists a solution of $g(x)=x$
in the open interval $(0,x^*)$ and the case $z<1/4$ would be resolved. We have $\dot g(x^*)=4(1-z)^2f(x^*)x^*=(2(1-z)x^*)^2$,
where we use $f(x^*)=x^*$. Then the inequality $\dot g(x)>1$ holds when $x^*>1/(2-2z)$ which after as simple algebra is reduced to
$z<1/4$. Thus when $z<1/4$ there are more than one fixed points of $T^2=T^2_{{\cal I},2}$.

Suppose now $z\geq 1/4$. Let again $x^*$ be  the unique solution to $f(x^*)=x^*$ and consider
 $\dot g(x)=4(1-z)^2f(x)x$. We claim
that when $z>1/4$, $\dot g(x)<1$ for all $x\in [0,1]$ and when $z=1/4$, $\dot g(x)<1$ for all $x\neq x^*$ and
$\dot g(x^*)=1$.
This immediately implies that $g(x)=x$ has at most one
solution, meaning it has exactly one since $g(x^*)=x^*$. Note $\dot g(0)=0<1$. We now prove that
$\dot g(1)=4(1-z)^2z<1$ and $\dot g(x)<1$ for all $x$ such that ${d^2g(x)\over dx^2}=0$, except for $x=x^*$
for which we will show that  $\dot g(x^*)=1$.

We start with $\dot g(1)=4(1-z)^2z\equiv \phi(z)$. Note $\phi(0)=\phi(1)=0<1$. Then the  maximum is
achieved at the points $z$  where $\dot\phi(z)=-8(1-z)z+4(1-z)^2=4(1-z)(1-3z)=0$. We have already
considered the case $z=1$. Otherwise $z=1/3$ for which $\phi(1/3)=16/27<1$. Thus $\dot g(1)\leq \sup_{z\in [0,1]}\phi(z)<1$.

Consider now points $x$ such that ${d^2g(x)\over dx^2}=4(1-z)^2[-2(1-z)x^2+1-(1-z)x^2]=0$, from which we obtain
$x=1/\sqrt{3-3z}$. Plugging this into
\[
\dot g(x)=4(1-z)^2(1-(1-z)x^2)x=4(1-z)^2{2\over 3}{1\over\sqrt{3-3z}}=({64(1-z)^3\over 27})^{1\over 2}.
\]
The last expression is smaller than  unity whenever $z>1/4$. When $z=1/4$ the expression is equal to the unity.
In this case, however, $x=1/\sqrt{3-3z}=2/3$ which is checked to be equal to $x^*={-1+\sqrt{5-4z} \over 2(1-z)}$ when $z=1/4$.
This concludes the proof of the proposition. \qed
\end{proof}

\section{Conclusions}\label{section:conclusions}
We have derived in this paper the limits of maximum weight independent sets and matchings in sparse random
graphs for some types of i.i.d weight distributions. Our method is based on a certain  \emph{local optimality} property
which states loosely that in  cases of certain distributions of the random weights, the optimal random
combinatorial structure under the consideration exhibits a long-range independence and, as a result, the value
which each node (edge) "contributes" to the optimal structure is almost completely determined by the constant
depth neighborhood of the node (edge). We certainly believe that such local optimality holds for many other
random combinatorial structures and it seems to be  an interesting property to study by itself, not to mention its
applications to studying random combinatorial structures.

\bibliographystyle{amsalpha}

\providecommand{\bysame}{\leavevmode\hbox
to3em{\hrulefill}\thinspace}
\providecommand{\MR}{\relax\ifhmode\unskip\space\fi MR }
\providecommand{\MRhref}[2]{%
  \href{http://www.ams.org/mathscinet-getitem?mr=#1}{#2}
} \providecommand{\href}[2]{#2}

\end{document}